\documentclass[12pt]{article}
\usepackage{amssymb}
\usepackage{amsmath}
\usepackage{color}
\usepackage{graphics}
\usepackage{epsfig}
\usepackage{hyperref}

\newtheorem{claim}{\bf \t}[part]

\newtheorem{Lemma}{Lemma}[part]

\newtheorem{Remark}{Remark}[part]
\newtheorem{Theorem}{Theorem}[part]

\numberwithin{Assumption}{section} \numberwithin{Corollary}{section}
\numberwithin{Definition}{section} \numberwithin{equation}{section}
\numberwithin{Example}{section} \numberwithin{Lemma}{section}
\numberwithin{Proposition}{section} \numberwithin{Remark}{section}
\numberwithin{Theorem}{section}

\def \Sum{\displaystyle\sum}

\def \ess{\mathrm{esssup}}

\def\t{\theta}

\def\text#1{{\rm #1}}
\begin{document}
\date{}
\title{\Large \bf Global existence and uniqueness of strong solutions to the 2D nonhomogeneous primitive  equations with density-dependent viscosity }

\author{\small \textbf{Quansen Jiu},$^{a}$
\thanks{
E-mail: jiuqs@mail.cnu.edu.cn.}\quad
\textbf{Lin Ma}$^{a}$,
\thanks{
E-mail: malin.cnu@foxmail.com.}\quad
and \textbf{Fengchao Wang}$^{b}$
\thanks{
E-mail: wfcwymm@163.com.}
}
 \maketitle
\small $^a$ School of Mathematical Sciences, Capital Normal University, Beijing
100048, P. R. China.

\small $^b$ College of Mathematics and Physics, Beijing University of Chemical Technology, Beijing 100029, P.R. China.

{\bf Abstract:}
This paper is concerned with an initial-boundary value problem of the  two-dimensional inhomogeneous  primitive  equations with density-dependent viscosity. The global well-posedness of strong solutions is established, provided the initial horizontal velocity is suitably small, that is, $\|\nabla u_{0}\|_{L^{2}}\leq \eta_{0}$ for suitably small $\eta_{0}>0$. The initial data may contain vacuum. The proof is based on the local well-posedness and the blow-up criterion proved in \cite{0}, which states that  if $T^{*}$ is the maximal existence time of the local strong solutions $(\rho,u,w,P)$ and $T^{*}<\infty$, then
\begin{equation*}
	\sup_{0\leq t<T^{*}}(\left\|\nabla \rho(t)\right\|_{L^{\infty}}+\left\|\nabla^{2}\rho(t)\right\|_{L^{2}}+\left\|\nabla u(t)\right\|_{L^{2}})=\infty.
\end{equation*}
To complete the proof, it is required to make an estimate on a key term
$\|\nabla u_{t}\|_{L_{t}^{1}L_{\Omega}^{2}}$. We prove that it is  bounded  and could be as small as desired
under certain smallness conditions, by making use of the regularity result of hydrostatic Stokes equations
and some careful time weighted estimates.


\section{Introduction } \setcounter{equation}{0}
\setcounter{Assumption}{0} \setcounter{Theorem}{0}
\setcounter{Proposition}{0} \setcounter{Corollary}{0}
\setcounter{Lemma}{0}

 \ \ \ \

In this paper, we consider the global existence and uniqueness of strong solutions to  the nonhomogeneous  primitive equations in $R^{2}$ , which read as
\begin{equation}\label{equ11}
\left\{
\begin{aligned}
&\rho_{t}+\partial_{x}(\rho u)+\partial_{y}(\rho w)=0,\\
&\partial_{t}(\rho u)+\partial_{x}(\rho u^{2})+\partial_{y}(\rho u w)+\partial_{x}P-\partial_{x}(\mu(\rho)\partial_{x}u)-\partial_{y}(\mu(\rho)\partial_{y}u)=0, \\
&\partial_{x}u+\partial_{y}w=0, \\
&\partial_{y}P=0,
\end{aligned}
\right.
\end{equation}
where the density $\rho$, the horizontal velocity $u$, the vertical velocity $w$ and the pressure $P$ are the unknown functions. The viscosity coefficient $\mu=\mu(\rho)$ is  a function of the density $\rho$ and assumed to satisfy
\begin{equation}\label{1}
\mu(\rho)\in C^{2}[0,\infty) \quad \text{and} \quad \mu(\rho)\geq\underline{\mu}>0 \quad\text{on}\quad [0,\infty).
\end{equation}

The primitive equations are based on the so called hydrostatic approximation, and they form a fundamental block in models for planetary oceanic and atmospheric dynamics (see, for example,
Pedlosky \cite{ped} and Vallis \cite{val}). As far as we know, it was essentially introduced by Richardson in \cite{rich}, to study the complicated atmospheric phenomena and for predicting the weather and possible climate changes.  Until 1990s, the mathematical studies of  primitive equations were initiated by  Lions, Temam and Wang in \cite{lions1992-1}-\cite{lions1995}, and they established the global existence of weak solutions to the system for the initial data $a\in L^{2}$ ; The uniqueness of weak solutions for the two-dimensional case was later proved by  Bresch, Kazhikjov and Lemoine \cite{bd} for initial data $a$ with $\partial_{z}a\in L^{2}$ and by  Kukavica, Pei, Rusin and Ziane \cite{Kuk2014} for the initial data which is  only continuous on the space variables. However, the uniqueness of these weak solutions for the three-dimensional case is still unclear. Concerning the strong solutions, the local existence results to the 2D primitive equations were proved by Guill-Gonzez, Masmoudi and Rodruez-Bellido in \cite{gui}, while the global existence and uniqueness of strong solutions for the 2D primitive equations was achieved by Bresch, Kazhikhov and Lemoine in \cite{bd} and Temam and Ziane in \cite{tem}. The global  existence of strong solutions for the 3D incompressible primitive equations has been known since the breakthrough work by Cao and Titi  \cite{cao1}, see also Kobelkov \cite{kob} and Kukavica and Ziane \cite{kuk7}-\cite{kuk8}.
In the partial viscosities and partial diffusivity case, one can see the works by Cao and Titi in \cite{cao2} and by Cao, Li and Titi in \cite{cao3}-\cite{cao6}. However,
the above-mentioned results are phrased within the $L^{2}$-setting. Based on the maximal regularity technique, Hieber and kashiwabara \cite{hie} and Giga, Gries, Hieber, Hussein and kashiwabara \cite{giga21}-\cite{giga20} gave the existing results in the
$L^{p}$ type spaces.

In recent years, the mathematical analysis of compressible primitive equations  has been attracting the attention of many mathematicians.
 Gatapov and Kazhikhov firstly studied in \cite{gag} and they proved the global existence of weak solutions to the 2D compressible primitive equations. The uniqueness of such  weak solutions was proved by Jiu, Li and Wang in \cite{jiu}.
Ersoy, Ngom and Sy \cite{ers2011} derived the 3D compressible primitive equations and obtained  the stability of weak solutions to the 3D compressible primitive equations with degenerate viscosity (see also Tang and Gao \cite{tang2015} for 2D case). Recently, Wang, Dou and Jiu \cite{wang} established the global existence of weak solutions to 3D compressible primitive equations with gravity field. The similar result is referred to Liu and Titi \cite{liu2019}. However, the uniqueness of these weak solutions for 3D case is still an open problem. Moreover, Liu and Titi \cite{liu2021} established local well-posedness results in the three-dimensional case. Meanwhile, they also considered the zero Mach number limit of compressible primitive equations in \cite{liu2020} and \cite{liu2023}.
 From the point of view of physics, the viscosity may depend on the density. Wang \cite{wang2021} established the local existence and uniqueness of  strong solutions to the 2D compressible primitive equations with density-dependent viscosity. Very recently, Wang and Jiu \cite{WJ2022} proved the global existence and uniqueness of strong solutions to the 2D compressible primitive equations with density-dependent viscosity.

In this paper, we consider the nonhomogeneous  primitive equations \eqref{equ11} in the domain $\Omega=M\times (0,1)$, with $M=(0,L)$, where $L>0$. We complement the system (\ref{equ11}) with boundary conditions and initial conditions as follows.
\begin{equation}\label{a2}
\begin{aligned}
 &\rho,~ u,~w,~ P\,~\text{are~periodic~ in~ the~direction}\,\, x,\\
 & u|_{y=0}=u|_{y=1}=0,\\
 &w|_{y=0}=w|_{y=1}=0,
\end{aligned}
\end{equation}
and
\begin{equation}\label{144}
\begin{aligned}
& u|_{t=0}=u_{0}(x,y),\\
& w|_{t=0}=w_{0}(x,y),\\
&\rho|_{t=0}=\rho_{0}(x,y).\\
\end{aligned}
\end{equation}

Throughout this paper, we define the following function spaces:
 \begin{align*}
 &V=\{v\in H^{1}(\Omega); \partial_{y}v=0\},\\
 & H^{1}_{0,\,per}(\Omega)=\{v\in H^1(\Omega); \partial_{x}\int_{0}^{1}v(x,y)dy=0, \,v|_{y=0,1}=0,\,v(0,y)=v(L,y)\}.
 \end{align*}

The nonhomogeneous primitive equations, as the intermediate equations between incompressible primitive equations and compressible primitive equations, were firstly studied by Jiu and Wang in $\cite{wang2}$, and they proposed the following compatibility condition:
\begin{equation*}\label{CC1}
  \mu\Delta u_{0}-\partial_{x}P_{0}=\rho^{\frac{1}{2}}_{0}g(x,z)\,\,\text{in}\,\, \Omega\,\,\,\, \text{and} \,\,\,\,\partial_{x}\int_{0}^{1}u_{0}dz=0,\quad \text{\,\,for\,\,some}\,\,(P_{0},g)\in V \times L^{2}
\end{equation*}
and established  the local existence and uniqueness of strong solutions to the 2D nonhomogeneous primitive equations with large data.
Furthermore, Wang, Jiu and Xu $\cite{wang1}$ obtianed an global existence result on strong solutions to the 2D nonhomogeneous primitive equations for initial data with small $H^{\frac{1}{2}}$-norm. Recently, for the case that the viscosity coefficient depends on the density, Jiu, Ma and Wang \cite{0}
established the local-in-time well-posedness of strong solutions to the 2D nonhomogeneous primitive equations \eqref{equ11}, which is listed as follows.
\begin{Theorem}\label{the-11}
Suppose that the initial data $(\rho_0,u_0)$ satisfies the regularity conditions
\begin{equation}\label{1711}
0\le \rho_{0}\in W^{2,2}(\Omega),~\nabla \rho_{0}\in L^{\infty},~u_{0}\in H^{1}_{0,per}(\Omega)\cap H^{2}(\Omega)
\end{equation} and the compatibility condition
\begin{equation}\label{188}	\begin{aligned}
		&-\rho_{0}u_{0}\partial_{x}u_{0}+\rho_{0}\int_{0}^{y}\partial_{x}u_{0}(s)ds\partial_{y}u_{0}-\partial_{x}P_{0}+\partial_{x}(\mu(\rho_{0})\partial_{x}u_{0})+\partial_{y}(\mu(\rho_{0})\partial_{y}u_{0}):=\rho_{0}V_{1},\\& \partial_{x}\int_{0}^{1}u_{0}(x,y)dy=0,
	\end{aligned}	
\end{equation}
for some $(P_{0},V_{1})\in V\times L^{2} $. Then there exists a small time $T_{*}\in(0,T)$, such that the initial boundary value problem $(\ref{equ11})$-$(\ref{144})$ has a unique strong solution $(\rho, u, w ,P)$ satisfying
\begin{equation}
\begin{aligned}
&\rho\geq 0,~\rho \in C([0,T_{*}];W^{1,\infty})\cap C([0,T_{*}];W^{2,2}),~\rho_{t}\in C([0,T_{*}];L^{q}),~ \sqrt{\rho}u_{t}\in L^{\infty}(0,T_{*};L^{2}),\\
&P\in C([0,T_{*}];H^{1})\cap L^{2}(0,T_{*};H^{2}),~
\nabla u\in C([0,T_{*}];H^{1})\cap L^{2}(0,T_{*};H^{2}),\\
&u_{t} \in L^{2}(0,T_{*};H^{1}_{0,per}),~
w\in C([0,T_{*}];H^{1})\cap L^{2}(0,T_{*};H^{2}),\\
\end{aligned}
\end{equation}
where $1\leqslant q<\infty$. Furthermore, if $T^{*}$  is the maximal existence time of the local strong solutions $(\rho,u,w,P)$, then
\begin{equation}\label{16}
\sup_{t\in[0,T^{*}]}(\left\|\nabla \rho(t)\right\|_{L^{\infty}}+\left\|\nabla^{2}\rho(t)\right\|_{L^{2}}+\left\|\nabla u(t)\right\|_{L^{2}})=\infty.
\end{equation}
\end{Theorem}

However, up to now, the global existence and uniqueness of strong solutions to the 2D nonhomogeneous primitive equations is open.
The main purpose of this paper is to establish the global well-posedness of strong solutions to the problem (\ref{equ11})-(\ref{144}) with small initial data in some sense. Our main result can be stated as

\begin{Theorem}\label{the-2}
Suppose that the initial data $(\rho_0,u_0,P_{0})$ satisfies the conditions in Theorem\ref{the-11}.
Then there exists a small positive constant $\eta_{0}$, depending on $\Omega,~\sup_{\Omega}\rho_{0}= \bar{\rho},~ \|\mu(\rho)\|_{C^{2}},~\|\nabla\rho_{0}\|_{L^{\infty}}$\\and $\|\nabla^{2}\rho_{0}\|_{L^{2}}$, 
such that if
$$\|\nabla u_{0}\|_{L^{2}}\leq \eta_{0},$$
then the initial boundary value problem $(\ref{equ11})$-$(\ref{144})$ admits a unique global strong solution\\ $(\rho,u,w,P)$ with
\begin{equation}
\begin{aligned}
&\rho\geq0,~\rho \in C(0,\infty;W^{1,\infty})\cap C(0,\infty;W^{2,2}),~\rho_{t}\in C(0,\infty;L^{q}),~ \sqrt{\rho}u_{t}\in L^{\infty}(0,\infty;L^{2}),\\
&P\in C(0,\infty;H^{1})\cap L^{2}(0,\infty;H^{2}),~\nabla u\in C(0,\infty;H^{1})\cap L^{2}(0,\infty;H^{2})\\
&u_{t} \in L^{2}(0,\infty;H^{1}_{0,per}),~ w\in C(0,\infty;H^{1})\cap L^{2}(0,\infty;H^{2}).
\end{aligned}
\end{equation}Moveover, the following decay property holds\begin{equation}\label{f}		\| \sqrt{\rho}u_{t}\|_{L^{2}}^{2}\leq C(\Omega,\|\mu\|_{C^{1}}, \bar{\rho},\|\nabla \rho_{0}\|_{L^{\infty}},\|\nabla^{2}\rho_{0}\|_{L^{2}})\|\nabla u_{0}\|_{L^{2}}^{2}t^{-2},	\end{equation}	and\begin{equation}\label{116}		\|\nabla u\|_{L^{2}}^{2}\leq C(\Omega,\|\mu\|_{C^{1}})(1+\frac{\overline{\rho}}{\|\rho_{0}\|_{L^{p}}})e^{-\sigma t},	\end{equation}	where $\sigma=\frac{\underline{\mu}}{\|\rho_{0}\|_{L^{p}}}$ for any $p>1$
.
\end{Theorem}


 Based on the result in \cite{0}, we apply the continuity method to establish the global well-posedeness of strong solutions to $(\ref{equ11})$-$(\ref{144})$. It was shown in \cite{0} that if $T^{*}$ is the maximal existence time of the local strong solutions $(\rho,u,w,P)$ and $T^{*}<T$, then
\begin{equation*}
	\sup_{0\leq t<T^{*}}(\left\|\nabla \rho(t)\right\|_{L^{\infty}}+\left\|\nabla^{2}\rho(t)\right\|_{L^{2}}+\left\|\nabla u(t)\right\|_{L^{2}})=\infty,
\end{equation*} which implies
that the proof of global existence relies on the estimates on the gradient of horizontal velocity and the one-order and two-order derivative of the density. We start with the
a priori hypothesises  $\|\nabla \rho\|_{L^{\infty}}+\|\nabla ^{2}\rho\|_{L^{2}}\leq8(\|\nabla \rho_{0}\|_{L^{\infty}}+\|\nabla^{2} \rho_{0}\|_{L^{2}})$ and $\|\nabla u\|_{L^{2}}^{2}\leq 4\frac{\bar{\mu}}{\underline {\mu}}\|\nabla u_{0}\|_{L^{2}}^{2},$ where $\bar{\mu}=\sup_{[0,\bar{\rho}]}\mu(\rho)$. Then we prove that $\|\nabla \rho\|_{L^{\infty}}+\|\nabla ^{2}\rho\|_{L^{2}}$ is less than $4(\|\nabla \rho_{0}\|_{L^{\infty}}+\|\nabla^{2} \rho_{0}\|_{L^{2}})$ and $\|\nabla u\|_{L^{2}}^{2}$ is less than $ (2ln2+1)\frac{\bar{\mu}}{\underline{\mu}}\|\nabla u_{0}\|_{L^{2}}^{2}$ on $[0,T]$ with the assumption $\|\nabla u_{0}\|_{L^{2}}$ is suitably small. In fact, by making use of the regularity
result of hydrostatic Stokes equations (see Lemma \ref{lem-33}), it is easy to obtain the desired bound of $\|\nabla u\|_{L^{2}}$ under the
assumption that $\|\nabla u_{0}\|_{L^{2}}$ is suitably small (see Lemma \ref{prop-1}). Then, the next key step is
to show that $\|\nabla \rho\|_{L^{\infty}}+\|\nabla ^{2}\rho\|_{L^{2}}$ is less than $4(\|\nabla \rho_{0}\|_{L^{\infty}}+\|\nabla^{2} \rho_{0}\|_{L^{2}})$. To do this, it
suffices to deal with a trouble term  $ \|\nabla u_{t}\|_{L_{t}^{1}L_{\Omega}^{2}}$ which   appears when estimating the  term $\|u\|_{L_{t}^{1}H_{\Omega}^{3}}$ (see Lemma \ref{the8}).
Indeed, based on the regularity result of hydrostatic Stokes equations
and some careful time weighted estimates, we find that the norm of
$\|\nabla u_{t}\|_{L_{t}^{1}L_{\Omega}^{2}}$ is bounded  and could be as small as desired
under certain smallness conditions. More precisely,
we split it into two terms \begin{equation}\label{3.93}\begin{aligned}	\int_{0}^{T_{*}}\|\nabla u_{t}\|_{L^{2}}dt&\leq(\int_{0}^{T^{0}}\|\nabla u_{t}\|_{L^{2}}^{2}dt)^{\frac{1}{2}}\sqrt{T^{0}}+(\int_{T^{0}}^{T_{*}}t^{2}\|\nabla u_{t}\|_{L^{2}}^{2}dt)^{\frac{1}{2}} (\int_{T^{0}}^{T_{*}}t^{-2}dt)^{\frac{1}{2}},
\end{aligned}
\end{equation}
where $T^{0}$ is a small enough number independent on the norm of $\|\nabla u_{0}\|_{L^{2}}.$ We remark that the local time $T_{*}$ might depend on $\|\nabla u_{0}\|_{L^{2}}$ (see \cite{0}). We can consider the local existence of the strong solutions under the additional assumption of $\|\nabla u_{0}\|_{L^{2}}\leq1$ such that the local time $T_{*}$ we obtain is independent on the norm of $\|\nabla u_{0}\|_{L^{2}}$ (see Lemma \ref{coro1}). Therefore, we can choose some small $T_{0}$ independent on the norm of $\|\nabla u_{0}\|_{L^{2}}$ such that the first term of (\ref{3.93}) is small enough. Meanwhile, the second term of (\ref{3.93}) can be controlled by the time weighted estimates ( see (\ref{3.34}) for more details).
Then, it follows that $\int_{0}^{T}\|u\|_{H^{3}}dt \leq \ln2$ under the assumption that $\|\nabla u_{0}\|_{L^{2}}$ is suitably small. 
 As a result, the control of $\|\nabla \rho\|_{L^{\infty}}+\|\nabla^{2}\rho \|_{L^{2}} $ and $\|\nabla u \|_{L^{2}}$ lead to uniform estimates for other higher order quantities, which guarantees the extension of local strong solutions.

 The rest of this paper is organized as follows. In Section 2, we recall some lemmas  and give an auxiliary lemma which is used later. In Section 3, we give an outline of the proof of the main Theorem.

 The following notations will be used:\\
 (1) $\|f\|_{L_{y}^{2}}=(\int_{0}^{1}|f|^{2}dy)^{\frac{1}{2}}$,\\
 (2) $\|f\|_{L_{y}^{2}L_{x}^{\infty}}=(\int_{0}^{1}\|f\|^{2}_{L_{x}^{\infty}}dy)^{\frac{1}{2}}$,
 $\|f\|_{L_{y}^{\infty}L_{x}^{\infty}}=\ess_{(x,y)\in \Omega} |f| $,\\
 (3)  For $ 1<p<\infty,\ \|f\|_{L^{p}}=\|f\|_{L_{y}^{p}L_{x}^{p}}=\|f\|_{L_{x}^{p}L_{y}^{p}}=(\int_{0}^{1}\int_{0}^{L}|f|^{p}dxdy)^{\frac{1}{p}}$,\\
 (4)$\int_{\Omega}=\int_{\Omega} dxdy$.
\section{Preliminaries}
\quad\quad
In this section, we recall some basic facts and elementary inequalities which will be used later. First of all, we present the well-known Gagliardo-Nirenberg interpolation inequality.
\begin{Lemma}\label{B1}
	(see $\cite{nir}$)(Gagliardo-Nirenberg interpolation inequality). For a function $u: \Omega \rightarrow \mathbb{R}$ defined on a bounded Lipschitz domain $\Omega \subset \mathbb{R}^{n}, \forall~  1 \leq q, r \leq \infty$, and a natural number m. Suppose also that a real number $\theta$ and a natural number $j$ are such that
	$$
	\frac{1}{p}=\frac{j}{n}+\left(\frac{1}{r}-\frac{m}{n}\right) \theta+\frac{1-\theta}{q}
	,$$
	and
	$$
	\frac{j}{m} \leq \theta \leq 1.
	$$
	Then, we have
	$$
	\left\|D^{j} u\right\|_{L^{p}} \leq C_{1}\left\|D^{m} u\right\|_{L^{r}}^{\theta}\|u\|_{L^{q}}^{1-\theta}+C_{2}\|u\|_{L^{s}},
	$$
	where s is a positive constant. The constants $C_{1}$ and $C_{2}$ depend on the domain $\Omega$ and $m, n, j, r, q, \theta$.
\end{Lemma}

The derivations of high-order estimates on the horizontal velocity $u$ rely on the following regularity results for the density-dependent hydrostatic Stokes equations.
\begin{Lemma}\label{B2}
	Assume that $\rho\in W^{2,2}$, $\nabla \rho \in L^{\infty}$ and $0\leq\rho\leq \bar{\rho}$, let $(u,P)\in H^{1}_{0,per}\times L^{2}$ be the unique weak solution to the problem
	\begin{equation}\label{l33}
		-\partial_{x}(\mu(\rho)\partial_{x}u)-\partial_{y}(\mu(\rho)\partial_{y}u)+\partial_{x}P=f, \quad \partial_{x}\int_{0}^{1}u(x,y)dy=0,\quad \int_{\Omega} P dxdy=0.
	\end{equation}
	Then we have the following regularity results:\\
	(1)If $f \in L^{2}$, then $(u, P) \in H^{2} \times H^{1}$ and
	$$
	\|u\|_{H^{2}}+\|P\|_{H^{1}} \leqslant C\left(1+\|\nabla \rho\|_{L^{\infty}}\right
	)\|f\|_{L^{2}}.
	$$
	(2)If $f \in L^{r}$, then $(u, P) \in W^{2,r} \times W^{1,r}$ and\begin{equation}
		\| u\|_{W^{2,r}}+\| P\|_{W^{1,r}}\leq C(1+\|\nabla\rho\|_{L^{\infty}})^{2}\|f\|_{L^{r}},
	\end{equation}
	where $r>2$.\\
	(3) If $f \in H^{1}$, then $(u, P) \in H^{3} \times H^{2}$ and
	\begin{equation}
		\|u\|_{H^{3}}+\|P\|_{H^{2}} \leq \widetilde{C}(1+\|\nabla\rho\|_{L^{ \infty}}+\|\nabla^{2}\rho\|_{L^{2}})^{3}\|f\|_{H^{1}},
	\end{equation}
	the constant $\widetilde{C}$ depend also on $\left\|\partial^{2} \mu / \partial x^{2}\right\|_{C}$.
\end{Lemma}
\begin{Lemma} \label{Lemma-2.3}
 		Assume that $u\in L^{1}([0,T],H^{3}(\Omega))$, let $\rho$ be the unique solution to the problem $$\partial_{t}\rho +\partial_{x}\rho u-\partial_{y}\rho\int_{0}^{y}\partial_{x}u(s)ds,~\rho|_{t=0}=\rho_{0}.$$ \\
 		If   $\nabla \rho_{0}\in L^{\infty}(\Omega)$ and $ \rho_{0}\in W^{2,2}(\Omega)$, then we have the following regularity results:
 		\begin{equation}\label{2.4}
 			\|\nabla \rho(t)\|_{L^{\infty}(\Omega)}\leq \|\nabla \rho_{0}\|_{L^{\infty}(\Omega)}\exp(\int_{0}^{t}\|u(s)\|_{H^{3}(\Omega)}ds),
 		\end{equation}
 	and
 		\begin{equation}\label{2.5}
 			\|\nabla ^{2}\rho(t)\|_{L^{2}(\Omega)}\leq (\|\nabla \rho_{0}\|_{L^{\infty}(\Omega)}+\|\nabla ^{2}\rho_{0}\|_{L^{2}(\Omega)})\exp\int_{0}^{t}\|u(s)\|_{H^{3}(\Omega)}ds.
 		\end{equation}
 \end{Lemma}

The proof of (\ref{2.4}) and (\ref{2.5})  is referred to \cite{0} and we omit it here.

As aforementioned, the global existence and uniqueness of strong solutions will be proved by combining the global a priori estimates in section 3  with the following local-in-time existence results which can follow as a corollary of Theorem \ref{the-11} in \cite{0}.
\begin{Lemma}\label{coro1}
 Suppose that the initial data $(\rho_0,u_0,P_{0
 })$ satisfies the regularity conditions in Theorem \ref{the-11}
and the additional condition $\|\nabla u_{0}\|_{L^{2}}^{2}\leq 1$. Then there exists a small time $T_{*}\in(0,T)$  such that the initial boundary value problem $(\ref{equ11})$-$(\ref{144})$ has a unique strong solution $(\rho, u, w, P)$
 satisfying	 	
\begin{equation*}
\begin{aligned}
&\quad \rho \in C([0,T_{*}];W^{1,\infty})\cap C([0,T_{*}];W^{2,2}),~\rho_{t}\in C([0,T_{*}];L^{q}),~ \sqrt{\rho}u_{t}\in L^{\infty}(0,T_{*};L^{2}),\\
&P\in C([0,T_{*}];H^{1})\cap L^{2}(0,T_{*};H^{2}),~	\nabla u\in C([0,T_{*}];H^{1})\cap L^{2}(0,T_{*};H^{2}),\\
&u_{t} \in L^{2}(0,T_{*};H^{1}_{0,per}),~w\in C([0,T_{*}];H^{1})\cap L^{2}(0,T_{*};H^{2}),\\
\end{aligned}
\end{equation*}	
where $1\leqslant q<\infty$ and $T_{*}$ is independent on the norm of $\|\nabla u_{0}\|_{L^{2}}$.
\end{Lemma}


\section{ Proof of Theorem \ref{the-2}}
\quad\quad 
This section is devoted to the proof of Theorem \ref{the-2}. Indeed, we follow the idea in \cite{0} and first prove the estimates for the case of a positive initial density and obtain the
uniform bounds which are independent of the lower bounds of the initial density. Then the vacuum case can be proved by applying standard regularizing techniques and compactness
arguments.
\subsection{A priori estimates}
\  \  \  \
 In this subsection, we establish several a priori estimates for strong solutions to the system $(\ref{equ11})$-$(\ref{144})$, which will play a key role in extending local strong solutions to the global ones. For simplicity, throughout this subsection, a constant $C$ will denote some positive one which may be dependent on $\Omega$ and $\|\mu\|_{C_{1}}$, but is independent of $\rho_{0}$ and $u_{0}$, and $\tilde{C}$ depends on $C$ and $\left\|\partial^{2} \mu / \partial x^{2}\right\|_{C}$. We remark that the value of the constant $C$ will be different from line to line.

We begin with the boundedness of density and the basic energy estimates,
which can be easily derived from the system $(\ref{equ11})$.

\begin{Lemma}\label{lemma-1}
	Assume that $(\rho, u, w, P)$ is the unique strong solution to the system $(\ref{equ11})$-$(\ref{144})$ on $\Omega\times(0,T)$ with the  initial data satisfying the assumptions in Theorem \ref{the-11}. Then, for every $(x,y,t) \in \Omega \times (0,T)$, it holds that
$$\|\rho\|_{L^{\infty}}=\|\rho_{0}\|_{L^{\infty}},$$
and
\begin{equation}\label{a}
\int_{\Omega} \rho \frac{1}{2}u^{2}(t)+\int_{0}^{t}\int_{\Omega} \mu(\rho)|\nabla u|^{2}ds\leq \frac{1}{2} \bar{\rho}\|u_{0}\|_{L^{2}}^{2},
\end{equation}
which yields
\begin{equation}\label{b}
\frac{1}{2}\int_{\Omega}\rho u^{2}(t)+\underline \mu\int_{0}^{t}\int_{\Omega}|\nabla u|^{2}ds\leq \frac{1}{2} \bar{\rho}\|u_{0}\|_{L^{2}}^{2}.
\end{equation}
\end{Lemma}

The next lemma is concerned with the $L^{2}$-norm of the gradient of horizontal velocity. Before sate and prove it, we denote
$$
M_{0}:=\|\nabla \rho_{0}\|_{L^{\infty}}+\|\nabla^{2} \rho_{0}\|_{L^{2}}.
$$
\begin{Lemma}\label{prop-1}
Assume that $(\rho, u , w, P)$ is the unique strong solution to the system $(\ref{equ11})$-$(\ref{144})$ on  $\Omega\times(0,T)$ with the  initial data satisfying the assumptions in Theorem \ref{the-11}, satisfying
\begin{equation}\label{311}
 \sup_{t\in[0,T]}\|\nabla \rho(t)\|_{L^{\infty}}+\sup_{t\in[0,T]}\|\nabla^{2}\rho(t)\|_{L^{2}}\leq 8M_{0},
 \end{equation}
 and
\begin{equation}\label{312}
\sup_{t\in[0,T]}\|\nabla u(t)\|_{L^{2}}^{2}\leq 4\frac{\bar{\mu}}{\underline{\mu}} \|\nabla u_{0}\|_{L^{2}}^{2}\leq 1.
\end{equation}
%
	 There exists a positive number $\gamma_{0}$, depending on $\Omega$, $\sup_{\Omega}\rho_{0}=\overline{\rho}$,  $\|\mu\|_{C^{2}}$, $\|\nabla\rho_{0}\|_{L^{\infty}}$ and $\| \nabla^{2} \rho_{0}\|_{L^{2}}$, such that, if
\begin{equation}\label{313}
	\|\nabla u_{0}\|_{L^{2}}\leq \gamma_{0},
\end{equation}
 then
\begin{equation}\label{314}
\frac{2}{\underline{\mu}}\int_{0}^{T}\int_{\Omega}\rho u^{2}_{t}dt+\sup_{t\in[0,T]}\int_{\Omega}|\nabla u|^{2}\leq(2ln2+1)\frac{\bar{\mu}}{\underline{\mu}}\|\nabla u_{0}\|_{L^{2}}^{2}.
\end{equation}
\end{Lemma}

Before proving Lemma \ref{prop-1}, we need the following lemma.
\begin{Lemma}\label{lem-33}
Assume that $(\rho, u, w, P)$ is the unique strong solution to the system $(\ref{equ11})$-$(\ref{144})$ on  $\Omega\times(0,T)$ with the  initial data satisfying the assumptions in Theorem \ref{the-11},, and satisfies

$$
 \sup_{t\in[0,T]}\|\nabla \rho(t)\|_{L^{\infty}}+\sup_{t\in[0,T]}\|\nabla^{2} \rho(t)\|_{L^{2}}\leq 8M_{0}.$$
Then it holds that
\begin{equation}\label{315}
\|\nabla u\|_{H^{1}}\leq C( \bar{\rho}^{1/2}M\| \sqrt\rho u_{t}\|_{L^{2}}+ \bar{\rho}M\|\nabla u\|^{2}_{L^{2}}+  \bar{\rho}^{\frac{3}{2}}M^{\frac{3}{2}}\|\nabla u\|_{L^{2}}^{\frac{5}{2}}+\bar{\rho}^{2} M^{2} \|\nabla u\|_{L^{2}}^{3}),
\end{equation}
where M =$1+8M_{0}$.
\end{Lemma}
{\bf Proof.}
 The horizontal momentum equations $\eqref{equ11}_{2}$ can be rewritten as follows,
\begin{equation}
-\partial_{x}(\mu(\rho)\partial_{x}u)
-\partial_{y}(\mu(\rho)\partial_{y}u)+\partial_{x}P=-\rho\partial_{t}u-\rho u\partial_{x}u+\rho\int_{0}^{y}\partial_{x}u(s)ds\partial_{y}u.
\end{equation}
It follows from Lemma $\ref{B1}$ and Lemma $\ref{B2}$ that
\begin{equation}\label{CC}
\begin{aligned}
\|\nabla u\|_{H^{1}}&\leq C  M(\| \rho u_{t}\|_{L^{2}}+\bar{\rho}\| u\|_{L^{6}}\|\partial_{x}u\|_{L^{3}}+\bar{\rho}\|\int_{0}^{y}\partial_{x} u(s)ds\|_{L_{y}^{2}L_{x}^{\infty}}\|\partial_{y}u\|_{L_{x}^{2}L_{y}^{\infty}})\\
&\leq CM[\bar{\rho}^{1/2}\| \sqrt{\rho} u_{t}\|_{L^{2}}+\bar{\rho}\|\nabla u\|_{L^{2}}(C_{1} \|\nabla^{2} u\|_{L^{2}}^{\frac{1}{3}}\|\nabla u\|_{L^{2}}^{\frac{2}{3}}+ C_{2}\|\nabla u\|_{L^{2}})\\
&\quad+\bar{\rho}\|\nabla u\|_{L^{2}}(C_{1} \|\nabla^{2} u\|_{L^{2}}^{\frac{1}{2}}\|\nabla u\|_{L^{2}}^{\frac{1}{2}}+ C_{2}\|\nabla u\|_{L^{2}})]\\
&\leq C( \bar{\rho}^{1/2} M\|\sqrt{\rho}  u_{t}\|_{L^{2}}+ \bar{\rho}M\|\nabla u\|^{2}_{L^{2}}+  \bar{\rho}M\|\nabla^{2} u\|_{L^{2}}^{\frac{1}{3}}\|\nabla u\|_{L^{2}}^{\frac{5}{3}}\\
&\quad+ M\bar{\rho}\|\nabla^{2} u\|_{L^{2}}^{\frac{1}{2}}\|\nabla u\|_{L^{2}}^{\frac{3}{2}}).
\end{aligned}
\end{equation}
Applying Young's inequality, we deduce that
\begin{equation}\label{318}
\|\nabla u\|_{H^{1}}\leq C( \bar{\rho}^{1/2}M\| \sqrt\rho u_{t}\|_{L^{2}}+ \bar{\rho}M\|\nabla u\|^{2}_{L^{2}}+  \bar{\rho}^{\frac{3}{2}}M^{\frac{3}{2}}\|\nabla u\|_{L^{2}}^{\frac{5}{2}}+\bar{\rho}^{2} M^{2} \|\nabla u\|_{L^{2}}^{3}).
\end{equation}$\hfill\Box$

Now we give the proof of Lemma \ref{prop-1}.

{\bf Proof of Lemma \ref{prop-1}.}
Multiplying the horizontal momentum equations $\eqref{equ11}_{2}$ by $u_{t}$ and integrating by part over $\Omega$ yields
\begin{equation}\label{C1}
\begin{aligned}
\int_{\Omega}\rho u^{2}_{t}+\frac{1}{2}\frac{d}{dt}\int_{\Omega}\mu(\rho)|\nabla u|^{2}&=-\int_{\Omega}\rho u \partial_{x}u u_{t}+\int_{\Omega}(\rho\int_{0}^{y}\partial_{x}u(s)ds \partial_{y}u u_{t})+\frac{1}{2}\int_{\Omega} \mu_{t}(\rho(t))|\nabla u|^{2} \\
&=-\int_{\Omega}\rho u \partial_{x}u u_{t}+\int_{\Omega}(\rho\int_{0}^{y}\partial_{x}u(s)ds\partial_{y}uu_{t})\\
&\quad+\frac{1}{2}\int_{\Omega}(\int_{0}^{y}\partial_{x}u(s)ds\partial_{y}\mu(\rho)-u(x)\partial_{x}\mu(\rho))|\nabla u||\nabla u|.
\end{aligned}
\end{equation}
Next, we will estimate the terms on the right-hand side in \eqref{C1} one by one.
For the first term, applying Young inequality and Lemma \ref{lem-33}, we have
\begin{equation}\label{L1}
\begin{aligned}
\int_{\Omega}\rho u_{t}u \partial_{x}u &\leq \|\sqrt\rho u_{t}\|_{L^{2}} \|\sqrt\rho \|_{L^{\infty}}\| u\|_{L^{6}}\|\partial_{x} u\|_{L^{3}}\\
&\leq \bar{\rho}^{\frac{1}{2}}\|\sqrt\rho u_{t}\|_{L^{2}}\|\nabla u\|_{L^{2}}(C_{1} \|\nabla^{2} u\|_{L^{2}}^{\frac{1}{3}}\|\nabla u\|_{L^{2}}^{\frac{2}{3}}+ C_{2}\|\nabla u\|_{L^{2}})\\
&\leq \epsilon \|\sqrt\rho u_{t}\|^{2}_{L^{2}}+\bar{\rho}\|\nabla u\|_{L^{2}}^{2}(C_{1} \|\nabla^{2} u\|_{L^{2}}^{\frac{2}{3}}\|\nabla u\|_{L^{2}}^{\frac{4}{3}}+ C_{2}\|\nabla u\|^{2}_{L^{2}})\\
&\leq \epsilon \|\sqrt\rho u_{t}\|^{2}_{L^{2}}+C(\bar{\rho}\|\nabla u\|_{L^{2}}^{\frac{10}{3}} \|\nabla^{2} u\|_{L^{2}}^{\frac{2}{3}}+\bar{\rho}\|\nabla u\|^{4}_{L^{2}})\\
&\leq \epsilon \|\sqrt\rho u_{t}\|^{2}_{L^{2}}+C( \bar{\rho}^{\frac{4}{3}} M^{\frac{2}{3}}\| \sqrt\rho u_{t}\|^{\frac{2}{3}}_{L^{2}}\|\nabla u\|^{\frac{10}{3}}+ M^{\frac{2}{3}}\bar{\rho}^{\frac{5}{3}}\|\nabla u\|^{\frac{14}{3}}_{L^{2}}\\
&\quad+ M\bar{\rho}^{2}\|\nabla u\|_{L^{2}}^{5}+ M^{\frac{4}{3}} \bar{\rho}^{\frac{7}{3}}\|\nabla u\|^{\frac{16}{3}}_{L^{2}}+\bar{\rho}\|\nabla u\|^{4}_{L^{2}})\\&\leq\epsilon \|\sqrt\rho u_{t}\|^{2}_{L^{2}}+C(\bar{\rho}^{2} M\|\nabla u\|^{5}+ M^{\frac{2}{3}}\bar{\rho}^{\frac{5}{3}}\|\nabla u\|^{\frac{14}{3}}_{L^{2}}\\
&\quad+ M^{\frac{4}{3}} \bar{\rho}^{\frac{7}{3}}\|\nabla u\|^{\frac{16}{3}}_{L^{2}}+\bar{\rho}\|\nabla u\|^{4}_{L^{2}}).
\end{aligned}
\end{equation}
Similarly, the second term can be bounded as
\begin{equation}\label{L2}
\begin{aligned}
\int_{\Omega}(\rho\int_{0}^{y}\partial_{x}u(s)ds\partial_{y}uu_{t})&
\leq \epsilon \|\sqrt\rho u_{t}\|^{2}_{L^{2}}+\bar{\rho}\|\nabla u\|_{L^{2}}^{2}(C_{1} \|\nabla^{2} u\|_{L^{2}}\|\nabla u\|_{L^{2}}+ C_{2}\|\nabla u\|^{2}_{L^{2}})\\
&\leq\epsilon \|\sqrt\rho u_{t}\|^{2}_{L^{2}}+C(\bar{\rho}\|\nabla u\|_{L^{2}}^{3} \|\nabla^{2} u\|_{L^{2}}+\bar{\rho}\|\nabla u\|^{4}_{L^{2}})\\
&\leq\epsilon \|\sqrt\rho u_{t}\|^{2}_{L^{2}}+C  ( \bar{\rho}^{3/2}M\| \sqrt\rho u_{t}\|_{L^{2}}\|\nabla u\|_{L^{2}}^{3}+ \bar{\rho}^{2}M\|\nabla u\|^{5}_{L^{2}}\\
&\quad+\bar{\rho}^{\frac{5}{2}}M^{\frac{3}{2}}\|\nabla u\|_{L^{2}}^{\frac{11}{2}}+\bar{\rho}^{3} M^{2} \|\nabla u\|_{L^{2}}^{6} +\bar{\rho}\|\nabla u\|^{4}_{L^{2}})\\
&\leq\epsilon \|\sqrt\rho u_{t}\|^{2}_{L^{2}}+C ( \bar{\rho}^{2}M\|\nabla u\|^{5}_{L^{2}}+  \bar{\rho}^{\frac{5}{2}}M^{\frac{3}{2}}\|\nabla u\|_{L^{2}}^{\frac{11}{2}}\\
&\quad+\bar{\rho}^{3} M^{2} \|\nabla u\|_{L^{2}}^{6} +\bar{\rho}\|\nabla u\|^{4}_{L^{2}}).
\end{aligned}
\end{equation}
For the last term, by direct estimates, we get
\begin{equation}\label{L3}
\begin{aligned}
&\int_{\Omega}[\int_{0}^{y}\partial_{x}u(s)ds\partial_{y}\mu(\rho)|\nabla u||\nabla u|-u(x)\partial_{x}\mu(\rho)|\nabla u||\nabla u|]\\
&\leq\|\int_{0}^{y}\partial_{x}u(s)ds\|_{L^{\infty}}\|\partial_{y}\mu(\rho)\|_{L^{\infty}}\|\nabla u\|_{L^{2}}\|\nabla u\|_{L^{2}}-\|u\|_{L^{\infty}}\|\partial_{x}\mu(\rho)\|_{L^{\infty}}\|\nabla u\|_{L^{2}}\|\nabla u\|_{L^{2}}\\
&\leq\|u\|_{H^{2}}\|\partial_{y}\mu(\rho)\|_{L^{\infty}}\|\nabla u\|_{L^{2}}^{2}\\
&\leq C( \bar{\rho}^{1/2}M\| \sqrt\rho u_{t}\|_{L^{2}}\|\nabla u\|_{L^{2}}^{2}+ \bar{\rho}M\|\nabla u\|^{4}_{L^{2}}+  \bar{\rho}^{\frac{3}{2}}M^{\frac{3}{2}}\|\nabla u\|_{L^{2}}^{\frac{9}{2}}+\bar{\rho}^{2} M^{2} \|\nabla u\|_{L^{2}}^{5}) M_{0}\\
&\leq\epsilon\| \sqrt\rho u_{t}\|_{L^{2}}^{2}+C(\bar{\rho} M^{2}M_{0}^{2}\|\nabla u\|_{L^{2}}^{4}+ \bar{\rho}MM_{0}\|\nabla u\|^{4}_{L^{2}}+ \bar{\rho}^{\frac{3}{2}}M^{\frac{3}{2}}M_{0}\|\nabla u\|_{L^{2}}^{\frac{9}{2}}+\bar{\rho}^{2} M^{2}M_{0} \|\nabla u\|_{L^{2}}^{5}).\\
\end{aligned}
\end{equation}
Substituting \eqref{L1}-\eqref{L3} into \eqref{C1}, we obtain
\begin{equation}\label{18}
	\begin{aligned}
		\int_{\Omega}\rho u^{2}_{t}+\frac{1}{2}\frac{d}{dt}\int_{\Omega}\mu(\rho)|\nabla u|^{2}
		&\leq 
C\mathcal{Q}\|\nabla u\|^{2}_{L^{2}},\\
	\end{aligned}
\end{equation}
where
\begin{equation}
\begin{aligned}
  \mathcal{Q}:=&\bar{\rho}\|\nabla u\|^{2}_{L^{2}}+\bar{\rho} M^{2}M_{0}^{2}\|\nabla u\|_{L^{2}}^{2}+ \bar{\rho}MM_{0}\|\nabla u\|^{2}_{L^{2}}+ M^{\frac{2}{3}}\bar{\rho}^{\frac{5}{3}}\|\nabla u\|^{\frac{8}{3}}_{L^{2}}\\
	&+ M^{\frac{4}{3}} \bar{\rho}^{\frac{7}{3}}\|\nabla u\|^{\frac{10}{3}}_{L^{2}}+ \bar{\rho}^{2}M\|\nabla u\|^{3}_{L^{2}}+\bar{\rho}^{2} M^{2}M_{0} \|\nabla u\|_{L^{2}}^{3}+ \bar{\rho}^{\frac{5}{2}}M^{\frac{3}{2}}\|\nabla u\|_{L^{2}}^{\frac{7}{2}}\\
 &+\bar{\rho}^{3} M^{2} \|\nabla u\|_{L^{2}}^{4} + \bar{\rho}^{\frac{3}{2}}M^{\frac{3}{2}}M_{0}\|\nabla u\|_{L^{2}}^{\frac{5}{2}}.
\end{aligned}
\end{equation}
Integrating with respect to time on $[0,t]$ and using the fact that  $\mu\geq\underline{\mu}>0$ gives
\begin{equation}\label{2011}
	\begin{aligned}
\frac{2}{\underline{\mu}}\int_{0}^{t}\int_{\Omega}\rho u^{2}_{t}ds+\int_{\Omega}|\nabla u|^{2}	
	&\leq\int_{0}^{t}\frac{C\mathcal{Q}}{\underline{\mu}}\|\nabla u(s)\|_{L^{2}}^{2}ds+\int_{\Omega}\frac{\mu(\rho_{0})}{\underline{\mu}}|\nabla u_{0}|^{2}\\&
	\leq\int_{0}^{t}
	\frac{C\mathcal{Q}}{\underline{\mu}}(\|\nabla u(s)\|_{L^{2}}^{2}+\frac{2}{\underline{\mu}}\int_{0}^{s}\int_{\Omega}\rho u_{t}^{2}(\tau)d\tau)ds+\int_{\Omega}\frac{\mu(\rho_{0})}{\underline{\mu}}|\nabla u_{0}|^{2}
	\end{aligned}
\end{equation}
Applying Gronwall inequality , we deduce that
\begin{equation}\label{20}
\begin{aligned}
	&\frac{2}{\underline{\mu}}\int_{0}^{T}\int_{\Omega}\rho u^{2}_{t}dt+\sup_{t\in[0,T]}\int_{\Omega}|\nabla u|^{2}\\
&\leq\frac{1}{\underline{\mu}}\int_{\Omega}\mu(\rho_{0})|\nabla u_{0}|^{2}+\frac{1}{\underline{\mu}}\int_{\Omega}\mu(\rho_{0})|\nabla u_{0}|^{2}\cdot\int_{0}^{T}\frac{C\mathcal{Q}}{\underline{\mu}}dt\cdot\exp(\int_{0}^{T}\frac{C\mathcal{Q}}{\underline{\mu}}dt).
	\end{aligned}
\end{equation}
According to Lemmas \ref{lemma-1},  ($\ref{312}$) and  Poincar\'e inequality, one has
 \begin{equation}\label{3.19}
\begin{aligned}
\int_{0}^{T}\frac{C\mathcal{Q}}{\underline{\mu}}dt &\leq C
\biggl\{\bar{\rho}+\bar{\rho} M^{2}M_{0}^{2}+ \bar{\rho}MM_{0}
+ M^{\frac{2}{3}}\bar{\rho}^{\frac{5}{3}}
+ M^{\frac{4}{3}} \bar{\rho}^{\frac{7}{3}}
+ \bar{\rho}^{2}M+\bar{\rho}^{2} M^{2}M_{0}\\
&\quad+\bar{\rho}^{\frac{5}{2}}M^{\frac{3}{2}}+\bar{\rho}^{3} M^{2} + \bar{\rho}^{\frac{3}{2}}M^{\frac{3}{2}}M_{0}
\biggr\}\frac{\bar {\rho}}{\underline{\mu}^{2}} \|\nabla u_{0}\|_{L^{2}}^{2}.
\end{aligned}
\end{equation}
Now we choose $\gamma_{0}$  small enough such that
\begin{equation*}\label{3133}
		\begin{aligned}
			C(\bar{\rho}+\bar{\rho} M^{4}+ \bar{\rho}M^{2 }		+ M^{\frac{2}{3}}\bar{\rho}^{\frac{5}{3}}	+ M^{\frac{4}{3}} \bar{\rho}^{\frac{7}{3}}	+ \bar{\rho}^{2}M+\bar{\rho}^{2} M^{3}	+  \bar{\rho}^{\frac{5}{2}}M^{\frac{3}{2}}+\bar{\rho}^{3} M^{2} + \bar{\rho}^{\frac{3}{2}}M^{\frac{5}{2}})\frac{\bar {\rho}}{\underline{\mu}^{2}} \gamma_{0}^{2} \leq \ln2.	\end{aligned}\end{equation*}
It  follows from  (\ref{20}), (\ref{3.19}) and (\ref{313}) that
\begin{equation}\label{}\begin{aligned}
\frac{2}{\underline{\mu}}\int_{0}^{T}\int_{\Omega}\rho u^{2}_{t}dt+\sup_{t\in[0,T]}\int_{\Omega}|\nabla u|^{2}&\leq\frac{\bar{\mu}}{\underline{\mu}}\|\nabla u_{0}\|^{2}_{L^{2}}+ \frac{\bar{\mu}}{\underline{\mu}}\|\nabla u_{0}\|^{2}\cdot ln2\cdot 2\\&=(2\ln2+1)\frac{\bar{\mu}}{\underline{\mu}}\|\nabla u_{0}\|_{L^{2}}.
	\end{aligned}
\end{equation}  $\hfill\Box$

Recall that $\|\nabla \rho\|_{L^{\infty}}+\|\nabla ^{2}\rho(t)\|_{L^{2}(\Omega)}\leq (\|\nabla \rho_{0}\|_{L^{\infty}(\Omega)}+\|\nabla ^{2}\rho_{0}\|_{L^{2}(\Omega)})\exp\int_{0}^{t}\|u(s)\|_{H^{3}(\Omega)}ds$. To prove that $\|\nabla \rho(t)\|_{L^{\infty
	}}+\|\nabla^{2} \rho(t)\|_{L^{2}}\leq 4M_{0}$, we still need the $H^{3}$-regularity of the velocity $u$. In fact, we need the following several lemmas to deal with $\|u\|_{H^{3}}$. First, we estiblish some time weighted estimates.
\begin{Lemma}\label{the5}
Assume that $(\rho, u, w, P)$ is the unique  strong solution to the system $(\ref{equ11})$-$(\ref{144})$ on $\Omega\times(0,T)$ with the initial data satisfying the assumptions in Theorem \ref{the-11} and (\ref{311})-(\ref{313}).  Then it holds that
\begin{equation}\label{331}
\frac{2}{\underline{\mu}}\int_{0}^{T}(t\int_{\Omega}\rho u_{t}^{2})dt+\sup_{t\in[0,T]}(t\int_{\Omega}|\nabla u|^{2})\leq(2\ln2+1)\frac{\bar{\mu}\bar{\rho}}{2\underline{\mu}^{2}}\|\nabla u_{0}\|_{L^{2}}^{2}.
\end{equation}

\end{Lemma}
{\bf Proof.}
Multiplying ($\ref{18}$) by $t$ and adding the term $\frac{1}{2}\int_{\Omega}\mu(\rho)|\nabla u|^{2}$, we have

\begin{equation}
\begin{aligned}
&t\int_{\Omega}\rho u_{t}^{2}+\frac{1}{2}\frac{d}{dt}(t\int_{\Omega}\mu(\rho)|\nabla u|^{2})
\leq tC\mathcal{Q}\|\nabla u\|_{L^{2 }}^{2}+\frac{1}{2}\int_{\Omega}\mu(\rho)|\nabla u|^{2}.\\
\end{aligned}
\end{equation}
Integrating with respect to time on $[0,t]$ and using the fact that  $\mu\geq\underline{\mu}>0$, we get
\begin{equation}
	\begin{aligned}
		&\frac{2}{\underline{\mu}}\int_{0}^{t}(s\int_{\Omega}\rho u_{t}^{2})ds +\sup_{t\in[0,T]}(t\int_{\Omega}|\nabla u|^{2})\\&\leq \int_{0}^{t}(
		\frac{C\mathcal{Q}}{\underline{\mu}}s\|\nabla u\|_{L^{2}}^{2})ds+\frac{1}{\underline{\mu}}\int_{0}^{t}(\int_{\Omega}\mu(\rho)|\nabla u|^{2})ds\\&\leq\int_{0}^{t}\frac{C\mathcal{Q}}{\underline{\mu}}(s\|\nabla u(s)\|_{L^{2}}^{2}+\frac{2}{\underline{\mu}}\int_{0}^{s}\tau\int_{\Omega}\rho u_{t}^{2}(\tau)d\tau)ds+\frac{1}{\underline{\mu}}\int_{0}^{t}(\int_{\Omega}\mu(\rho)|\nabla u|^{2})ds.
	\end{aligned}
\end{equation}
By Gronwall inequality, assumption (\ref{313}), Poincar\'e inequality and Lemma \ref{lemma-1} , we deduce
\begin{equation}
\begin{aligned}
&\frac{2}{\underline{\mu}}\int_{0}^{T}(t\int_{\Omega}\rho u_{t}^{2})dt+\sup_{t\in[0,T]}(t\int_{\Omega}|\nabla u|^{2})\\
&\leq\frac{1}{\underline{\mu}}\int_{0}^{T}(\int_{\Omega}\mu(\rho)|\nabla u|^{2})dt+\frac{1}{\underline{\mu}}\int_{0}^{T}(\int_{\Omega}\mu(\rho)|\nabla u|^{2})dt\cdot\int_{0}^{T}\frac{C\mathcal{Q}}{\underline
	{\mu}}dt\cdot\exp(\int_{0}^{T}\frac{C\mathcal{Q}}{\underline
{\mu}}dt)\\&\leq(2\ln2+1)\frac{\bar{\mu}\bar{\rho}}{2\underline{\mu}^{2}}\|\nabla u_{0}\|_{L^{2}}^{2}.
\end{aligned}
\end{equation}$\hfill\Box$
\begin{Lemma}\label{the6}
Assume that $(\rho, u, w, P)$ is the unique strong solution to the system $(\ref{equ11})$-$(\ref{144})$ on $\Omega\times(0,T)$ with the initial data satisfying the assumptions in Theorem \ref{the-11} and (\ref{311})-(\ref{313}). Then it holds that
\begin{equation}\label{333}
\begin{aligned}
&\sup_{t\in[0,T]}t\int_{\Omega} \rho u_{t}^{2}+2\underline{\mu}\int_{0}^{T}(t\int_{\Omega}|\nabla u_{t}|^{2})dt\leq C \Theta_{1}\|\nabla u_{0}\|_{L^{2}}^{2}\cdot\exp\big[ C\Theta_{2}\big],
\end{aligned}
\end{equation}
and
\begin{equation}\label{3.34}
\begin{aligned}
&\sup_{t\in[0,T]}t^{2}\int_{\Omega} \rho u_{t}^{2}+2\underline{\mu}\int_{0}^{T}(t^{2}\int|\nabla u_{t}|^{2})dt	\leq
C \bar{\rho}\Theta_{1}\|\nabla u_{0}\|_{L^{2}}^{2}\cdot\exp\big[ C\Theta_{2}\big],
\end{aligned}
\end{equation}
where $\Theta_{1}=1+{M}^{4}\bar{\rho}^{4}+
\bar{\rho}^{5}  M^{5}+ {M}^{6}\bar{\rho}^{6}$ and $\Theta_{2}=C( {M}^{4}\bar{\rho}^{2}+
\bar{\rho}^{2} + \bar{\rho}^{4})$.

\end{Lemma}
{\bf Proof.}
Taking t-derivative of the horizontal momentum equations $\eqref{equ11}_{2}$, we deduce that
\begin{equation}\label{335}
\begin{aligned}
&\rho u_{tt}+\rho u \partial_{x}u_{t}-\rho\int_{0}^{y}\partial_{x}u(s)ds \partial_{y}u_{t}-\partial_{x}(\mu(\rho)\partial_{x}u_{t})-\partial_{y}(\mu(\rho)\partial_{y}u_{t})+\partial_{x}P_{t}\\
&=-\rho_{t}u_{t}-(\rho u)_{t}\partial_{x}u+(\rho\int_{0}^{y}\partial_{x}u(s)ds)_{t}\partial_{y}u+\partial_{x}(\mu(\rho)_{t}\partial_{x}u)+\partial_{y}(\mu(\rho)_{t}\partial_{y}u).
\end{aligned}
\end{equation}
Multiplying  (\ref{335}) by $tu_{t}$ and integrating (by parts) over $\Omega$ yields
\begin{equation}\label{336}
\begin{aligned}
&\frac{t}{2}\frac{d}{dt}\int_{\Omega}\rho u_{t}^{2}+t\int_{\Omega}\mu(\rho)|\nabla u_{t}|^{2}\\
&=-\frac{t}{2}\int_{\Omega}\rho_{t}u_{t}^{2}-t\int_{\Omega}(\rho u)_{t}\partial_{x}u u_{t}+t\int_{\Omega}(\rho\int_{0}^{y}\partial_{x}u(s)ds)_{t}\partial_{y}u u_{t}-\frac{t}{2}\int_{\Omega}(\mu(\rho))_{t}\frac{d|\nabla u|^{2}}{dt}.
\end{aligned}
\end{equation}
Noting the fact that $0<\underline{\mu}\leq \mu$, one has
\begin{equation}\label{342}
\begin{aligned}
&\frac{t}{2\underline{\mu}}\frac{d}{dt}\int_{\Omega}\rho u_{t}^{2}+t\int_{\Omega}|\nabla u_{t}|^{2}=-\frac{t}{2\underline{\mu}}\int_{\Omega}\rho_{t}u_{t}^{2}-\frac{t}{\underline{\mu}}\int_{\Omega}(\rho u)_{t}\partial_{x}u u_{t}\\
&+\frac{t}{\underline{\mu}}\int_{\Omega}(\rho\int_{0}^{y}\partial_{x}u(s)ds)_{t}\partial_{y}u u_{t}-\frac{t}{2\underline{\mu}}\int_{\Omega}(\mu(\rho))_{t}\frac{d|\nabla u|^{2}}{dt}=\sum_{i=1}^{4}I_{i}.
\end{aligned}
\end{equation}
 Using $(\ref{equ11})_1$ and Poincar\'e inequality, we estimate the term $I_{1}$ as
\begin{equation}
\begin{aligned}
I_{1}&=-\frac{t}{2\underline{\mu}}\int_{\Omega}\rho_{t}u_{t}^{2}\\&=\frac{t}{2\underline{\mu}}\int_{\Omega}[\partial_{x}(\rho u)-\partial_{y}(\rho\int_{0}^{y}\partial_{x}u(s)ds)]u_{t}u_{t}\\
&=-\frac{t}{\underline{\mu}}\int_{\Omega} \rho u \partial_{x}u_{t}u_{t}+\frac{t}{\underline{\mu}}\int_{\Omega} \rho\int_{0}^{y}\partial_{x}u(s)ds\partial_{y}u_{t}u_{t}\\
&\leq \frac{t}{\underline{\mu}}\|\rho u_{t}\|_{L^{2}_{y}L^{\infty}_{x}}\|\partial_{x}u_{t}\|_{L^{2}}\|\nabla u\|_{L^{2}}+\frac{t}{\underline{\mu}}\|\rho u_{t}\|_{L^{2}_{y}L^{\infty}_{x}}\|\int_{0}^{y}\partial_{x}u(s)ds\|_{L^{2}_{x}L^{\infty}_{y}}\|\partial_{y}u_{t}\|_{L^{2}}\\
&\leq C\frac{t}{\underline{\mu}}(\|\partial_{x}\rho u_{t}+\rho \partial_{x}u_{t}\|_{L^{2}}^{\frac{1}{2}}\|\rho u_{t}\|_{L^{2}}^{\frac{1}{2}}+\|\rho u_{t}\|_{L^{2}})\|\nabla u_{t}\|_{L^{2}}\|\nabla u\|_{L^{2}}\\
&\leq C\frac{t}{\underline{\mu}}[(M_{0}\|u_{t}\|_{L^{2}}^{\frac{1}{2}}+\bar{\rho}^{\frac{1}{2}}\|\partial_{x}u_{t}\|_{L^{2}}^{\frac{1}{2}})\|\rho u_{t}\|^{\frac{1}{2}}_{L^{2}}+\|\rho u_{t}\|_{L^{2}}]\|\nabla u_{t}\|_{L^{2}}\|\nabla u\|_{L^{2}}\\
&\leq C(\frac{t}{\underline{\mu}}\bar{\rho}^{\frac{1}{4}}M_{0}\|\nabla u_{t}\|_{L^{2}}^{\frac{3}{2}}\|\nabla u\|_{L^{2}}\|\sqrt{\rho} u_{t}\|_{L^{2}}^{\frac{1}{2}}+\frac{t}{\underline{\mu}}\bar{\rho}^{\frac{3}{4}}\|\nabla u_{t}\|_{L^{2}}^{\frac{3}{2}}\|\nabla u\|_{L^{2}}\|\sqrt{\rho} u_{t}\|_{L^{2}}^{\frac{1}{2}}\\
&\quad+\frac{t}{\underline{\mu}}\bar{\rho}^{\frac{1}{2}}\|\nabla u_{t}\|_{L^{2}}\|\nabla u\|_{L^{2}}\|\sqrt{\rho} u_{t}\|_{L^{2}}).\\
\end{aligned}
\end{equation}
Applying Young's inequality, we deduce that
\begin{equation}\label{3399}
\begin{aligned}
I_{1}&\leq \epsilon t\|\nabla u_{t}\|^{2}+Ct\bar{\rho} M_{0}^{4}\|\sqrt{\rho} u_{t}\|^{2}\|\nabla u\|^{4}+Ct\bar{\rho}^{3}\|\sqrt{\rho} u_{t}\|^{2}\|\nabla u\|^{4}\\
&\quad+Ct\bar{\rho}\|\sqrt{\rho} u_{t}\|^{2}\|\nabla u\|^{2},
\end{aligned}
\end{equation}
where $\epsilon$ is some small number.

To estimate $I_{2}$, taking the mass equation  $(\ref{equ11})_1$ into consideration, we get
\begin{equation*}
\begin{aligned}
I_{2}&=-\frac{t}{\underline{\mu}}\int_{\Omega}(\rho u)_{t}\partial_{x}u u_{t}=-\frac{t}{\underline{\mu}}\int_{\Omega}\rho_{t} u\partial_{x}u u_{t}-\frac{t}{\underline{\mu}}\int_{\Omega}\rho u_{t}\partial_{x}u u_{t}\\
&=\frac{t}{\underline{\mu}}\int_{\Omega}[\partial_{x}\rho u-\partial_{y}\rho\int_{0}^{y}\partial_{x}u(s)ds]u \partial_{x}u u_{t}-\frac{t}{\underline{\mu}}\int_{\Omega}\rho u_{t}\partial_{x}u u_{t}
= I_{21}+I_{22}.
\end{aligned}
\end{equation*}
It follows from Sobolev embedding inequality, Gagliardo-Nirenberg inequality, and Lemma \ref{lem-33} that
\begin{equation}
\begin{aligned}
I_{21}&=\frac{t}{\underline{\mu}}\int_{\Omega}[\partial_{x}\rho u u \partial_{x}u u_{t}-\partial_{y}\rho\int_{0}^{y}\partial_{x}u(s)ds u \partial_{x}u u_{t}]\\\
&\leq\frac{t}{\underline{\mu}}\|\partial_{x}\rho\|_{L^{\infty}}\|u\|_{L^{\infty}}\|u\|_{L^{2}_{x}L^{\infty}_{y}}\|\partial_{x}u\|_{L^{2}}\|u_{t}\|_{L^{2}_{y}L^{\infty}_{x}}\\
&\quad+\frac{t}{\underline{\mu}} \|\partial_{y}\rho\|_{L^{\infty}}\|u\|_{L^{\infty}}\|\int_{0}^{y}\partial_{x}u(s)ds\|_{L^{2}_{x}L^{\infty}_{y}}\|\partial_{x}u\|_{L^{2}}\|u_{t}\|_{L^{2}_{y}L^{\infty}_{x}}\\
&\leq\frac{t}{\underline{\mu}}M_{0}\|u\|_{H^{2}}\|\nabla u\|_{L^{2}}^{2}\|\nabla u_{t}\|_{L^{2}}\\
&\leq C\frac{t}{\underline{\mu}}M_{0}(\bar{\rho}^{1/2}M\| \sqrt\rho u_{t}\|_{L^{2}}+ \bar{\rho}M\|\nabla u\|^{2}_{L^{2}}+  \bar{\rho}^{\frac{3}{2}}M^{\frac{3}{2}}\|\nabla u\|_{L^{2}}^{\frac{5}{2}}+ \bar{\rho}^{2} M^{2} \|\nabla u\|_{L^{2}}^{3} )\|\nabla u\|^{2}_{L^{2}}\|\nabla u_{t}\|_{L^{2}}\\
&\leq \epsilon t\|\nabla u_{t}\|^{2}+Ct(M_{0}^{2}\bar{\rho}  M^{2}\| \sqrt\rho u_{t}\|_{L^{2}}^{2}\|\nabla u\|^{4}+M_{0}^{2}M^{2}\bar{\rho}^{2}\|\nabla u\|^{8}_{L^{2}}+M_{0}^{2} M^{3}\bar{\rho}^{3}\|\nabla u\|^{9}_{L^{2}}\\
&\quad+M_{0}^{2} M^{4}\bar{\rho}^{4}\|\nabla u\|^{10}_{L^{2}}),
\end{aligned}
\end{equation}
where $\epsilon$ is some small number. Similarly,  it holds that
\begin{equation}
\begin{aligned}
I_{22}&=-\frac{t}{\underline{\mu}}\int_{\Omega}\rho u_{t}\partial_{x}u u_{t}\\&\leq \frac{t}{\underline{\mu}}\|\rho u_{t}\|_{L^{2}_{y}L^{\infty}_{x}}\|\nabla u\|_{L^{2}}\|\nabla u_{t}\|_{L^{2}}\\
&\leq \epsilon t\|\nabla u_{t}\|^{2}+Ct \bar{\rho}M_{0}^{4}\|\sqrt{\rho} u_{t}\|^{2}\|\nabla u\|^{4}+Ct\bar{\rho}^{3}\|\sqrt{\rho} u_{t}\|^{2}\|\nabla u\|^{4}
+Ct\bar{\rho}\|\sqrt{\rho} u_{t}\|^{2}\|\nabla u\|^{2},
\end{aligned}
\end{equation}
where $\epsilon$ is some small number. 

To bound the term $I_{3}$, we rewrite it as
\begin{equation*}
\begin{aligned}
I_{3}&=\frac{t}{\underline{\mu}}\int_{\Omega}(\rho\int_{0}^{y}\partial_{x}u(s)ds)_{t}\partial_{y}u u_{t}\\
&=\frac{t}{\underline{\mu}}\int_{\Omega} [\rho_{t}\int_{0}^{y}\partial_{x}u(s)ds\partial_{y}u u_{t}]+\frac{t}{\underline{\mu}}\int_{\Omega} [\rho\int_{0}^{y}\partial_{x}u_{t}(s)ds\partial_{y}u u_{t}]\\
&=-\frac{t}{\underline{\mu}}\int_{\Omega}[(\partial_{x}\rho u-\partial_{y}\rho\int_{0}^{y}\partial_{x}u(s)ds)\int_{0}^{y}\partial_{x}u(s)ds\partial_{y}u u_{t}]+\frac{t}{\underline{\mu}}\int_{\Omega} [\rho\int_{0}^{y}\partial_{x}u_{t}(s)ds\partial_{y}u u_{t}]\\
&=I_{31}+I_{32}.
\end{aligned}
\end{equation*}
Owing to Young's inequality and Sobolev embedding inequality, we deduce that
\begin{equation}
\begin{aligned}
I_{31}&\leq\frac{t}{\underline{\mu}}\|\partial_{x}\rho\|_{L^{\infty}}\|u\|_{L^{\infty}}\|\int_{0}^{y}\partial_{x}u(s)ds\|_{L^{\infty}_{y}L^{2}_{x}}
\|\partial_{y}u\|_{L^{2}}\| u_{t}\|_{L^{\infty}_{x}L^{2}_{y}}\\
&\quad+\frac{t}{\underline{\mu}}\|\partial_{y}\rho\|_{L^{\infty}}\|\int_{0}^{y}\partial_{x}u(s)ds\|_{L^{\infty}}\|\int_{0}^{y}\partial_{x}u(s)ds\|_{L^{\infty}_{y}L^{2}_{x}}
\|\partial_{y}u\|_{L^{2}}\| u_{t}\|_{L^{\infty}_{x}L^{2}_{y}}\\
&\leq \frac{t}{\underline{\mu}}M_{0}\|u\|_{H^{2}}\|\nabla u\|_{L^{2}}^{2}\|\nabla u_{t}\|_{L^{2}}\\
&\leq \epsilon t\|\nabla u_{t}\|_{L^{2}}^{2}+Ct
(M_{0}^{2}\bar{\rho}  M^{2}\| \sqrt\rho u_{t}\|_{L^{2}}^{2}\|\nabla u\|_{L^{2}}^{4}+M_{0}^{2}M^{2}\bar{\rho}^{2}\|\nabla u\|^{8}_{L^{2}}\\
&\quad+M_{0}^{2} M^{3}\bar{\rho}^{3}\|\nabla u\|^{9}_{L^{2}}+M_{0}^{2} M^{4}\bar{\rho}^{4}\|\nabla u\|^{10}_{L^{2}}),
\end{aligned}
\end{equation}
where $\epsilon$ is some small number. A similar argument to $I_{32}$ yields
\begin{equation}
\begin{aligned}
I_{32}&=\frac{t}{\underline{\mu}}\int_{\Omega}[\rho\int_{0}^{y}\partial_{x}u_{t}(s)ds\partial_{y}u u_{t}] \\
&\leq\frac{t}{\underline{\mu}}\|\rho u_{t}\|_{L^{2}_{y}L^{\infty}_{x}}\|\int_{0}^{y}\partial_{x}u_{t}(s)ds\|_{L^{2}_{x}L^{\infty}_{y}}\|\partial_{y}u\|_{L^{2}}\\
&\leq \frac{t}{\underline{\mu}}\|\rho u_{t}\|_{L^{2}_{y}L^{\infty}_{x}}\|\nabla u\|_{L^{2}}\|\nabla u_{t}\|_{L^{2}}\\
&\leq \epsilon t\|\nabla u_{t}\|^{2}+Ct \bar{\rho}M_{0}^{4}\|\sqrt{\rho} u_{t}\|^{2}\|\nabla u\|^{4}+Ct\bar{\rho}^{3}\|\sqrt{\rho} u_{t}\|^{2}\|\nabla u\|^{4}
+Ct\bar{\rho}\|\sqrt{\rho} u_{t}\|^{2}\|\nabla u\|^{2},
\end{aligned}
\end{equation}
where $\epsilon$ is some small number.

Concerning the last term $I_{4}$, by direct estimates, one has
\begin{equation}\label{3466}
\begin{aligned}
I_{4}&=-\frac{t}{\underline{\mu}}\int_{\Omega}\mu(\rho)_{t}\nabla u\nabla u_{t}\\
&=\frac{t}{\underline{\mu}}\int_{\Omega}[\partial_{x}\mu(\rho)u-\partial_{y}\mu(\rho)\int_{0}^{y}\partial_{x}u(s)ds]\nabla u\nabla u_{t}\\
&\leq\frac{t}{\underline{\mu}}\|\nabla\mu(\rho)\|_{L^{\infty}}\|u\|_{H^{2}}\|\nabla u\|_{L^{2}}\|\nabla u_{t}\|_{L^{2}}\\&\leq C\frac{t}{\underline{\mu}}M_{0}( \bar{\rho}^{1/2}M\| \sqrt\rho u_{t}\|_{L^{2}}+ \bar{\rho}M\|\nabla u\|^{2}_{L^{2}}+  \bar{\rho}^{\frac{3}{2}}M^{\frac{3}{2}}\|\nabla u\|_{L^{2}}^{\frac{5}{2}}+\bar{\rho}^{2} M^{2} \|\nabla u\|_{L^{2}}^{3})\|\nabla u\|_{L^{2}}\|\nabla u_{t}\|_{L^{2}}\\
&\leq \epsilon t\|\nabla u_{t}\|^{2}+Ct(M_{0}^{2}\bar{\rho}  M^{2}\| \sqrt\rho u_{t}\|_{L^{2}}^{2}\|\nabla u\|_{L^{2}}^{2}+M_{0}^{2}M^{2}\bar{\rho}^{2}\|\nabla u\|^{6}_{L^{2}}+M_{0}^{2} M^{3}\bar{\rho}^{3}\|\nabla u\|^{7}_{L^{2}}\\
&\quad+M_{0}^{2} M^{4}\bar{\rho}^{4}\|\nabla u\|^{8}_{L^{2}}),
\end{aligned}
\end{equation}
where $\epsilon$ is some small number.

Substituting all the above estimates $(\ref{3399})$-$(\ref{3466})$ into $(\ref{342})$ and choosing $\epsilon$ small enough, we deduce that
\begin{equation}\label{347}
\begin{aligned}
&\frac{t}{2\underline{\mu}}\frac{d}{dt}\int_{\Omega}\rho u_{t}^{2}+t\int_{\Omega}|\nabla u_{t}|^{2}\\
&\leq Ct( \bar{\rho}M_{0}^{4}\|\sqrt{\rho} u_{t}\|_{L^{2}}^{2}\|\nabla u\|^{4}_{L^{2}}+\bar{\rho}^{3}\|\sqrt{\rho} u_{t}\|_{L^{2}}^{2}\|\nabla u\|_{L^{2}}^{4}+\bar{\rho}\|\sqrt{\rho} u_{t}\|_{L^{2}}^{2}\|\nabla u\|_{L^{2}}^{2}\\
&\quad+M_{0}^{2}\bar{\rho}  M^{2}\| \sqrt\rho u_{t}\|_{L^{2}}^{2}\|\nabla u\|_{L^{2}}^{2}+M_{0}^{2}M^{2}\bar{\rho}^{2}\|\nabla u\|^{6}_{L^{2}}+M_{0}^{2} M^{3}\bar{\rho}^{3}\|\nabla u\|^{7}_{L^{2}}+M_{0}^{2} M^{4}\bar{\rho}^{4}\|\nabla u\|^{8}_{L^{2}}\\
&\quad+M_{0}^{2}\bar{\rho}  M^{2}\| \sqrt\rho u_{t}\|_{L^{2}}^{2}\|\nabla u\|^{4}+M_{0}^{2}M^{2}\bar{\rho}^{2}\|\nabla u\|^{8}_{L^{2}}+M_{0}^{2} M^{3}\bar{\rho}^{3}\|\nabla u\|^{9}_{L^{2}}+M_{0}^{2} M^{4}\bar{\rho}^{4}\|\nabla u\|^{10}_{L^{2}}).
\end{aligned}
\end{equation}
Adding the term $\frac{1}{2\underline{\mu}}\int_{\Omega}\rho u_{t}^{2}$ into $(\ref{347})$, one has
\begin{equation}\label{346}
\begin{aligned}
&\frac{1}{2\underline{\mu}}\frac{d}{dt}(t\int_{\Omega} \rho u_{t}^{2})+t\int_{\Omega}|\nabla u_{t}|^{2}\\
&\leq C t\| \sqrt\rho u_{t}\|_{L^{2}}^{2}( \bar{\rho}M_{0}^{4}\|\nabla u\|_{L^{2}}^{4}+\bar{\rho}^{3}\|\nabla u\|_{L^{2}}^{4}
+\bar{\rho}\|\nabla u\|_{L^{2}}^{2}+M_{0}^{2}\bar{\rho}  M^{2}\|\nabla u\|_{L^{2}}^{2}+
M_{0}^{2}\bar{\rho}  M^{2}\|\nabla u\|_{L^{2}}^{4})\\
&\quad+Ct(M_{0}^{2}M^{2}\bar{\rho}^{2}\|\nabla u\|^{6}_{L^{2}}+M_{0}^{2} M^{3}\bar{\rho}^{3}\|\nabla u\|^{7}_{L^{2}}+M_{0}^{2} M^{4}\bar{\rho}^{4}\|\nabla u\|^{8}_{L^{2}}
+M_{0}^{2}M^{2}\bar{\rho}^{2}\|\nabla u\|^{8}_{L^{2}}\\
&\quad+M_{0}^{2} M^{3}\bar{\rho}^{3}\|\nabla u\|^{9}_{L^{2}}+M_{0}^{2} M^{4}\bar{\rho}^{4}\|\nabla u\|^{10}_{L^{2}})
+\frac{1}{2\underline{\mu}}\|\sqrt\rho u_{t}\|_{L^{2}}^{2}.
\end{aligned}
\end{equation}
Integrating with respect to time on $[0,t]$ and using the fact that  $\mu\geq\underline{\mu}>0$, we get
\begin{equation*}
	\begin{aligned}
	&t\int_{\Omega}\rho u_{t}^{2}+2\underline{\mu}\int_{0}^{t}s\int_{\Omega}|\nabla u_{t}|^{2}ds\\
&\leq\int_{0}^{t}Cs\|\sqrt{\rho}u_{t}\|_{L^{2}}^{2}( \bar{\rho}M_{0}^{4}\|\nabla u\|_{L^{2}}^{4}+\bar{\rho}^{3}\|\nabla u\|_{L^{2}}^{4}
	+\bar{\rho}\|\nabla u\|_{L^{2}}^{2}+M_{0}^{2}\bar{\rho}  M^{2}\|\nabla u\|_{L^{2}}^{2}\\&\quad+
	M_{0}^{2}\bar{\rho}  M^{2}\|\nabla u\|_{L^{2}}^{4})ds+\int_{0}^{t}Cs(M_{0}^{2}M^{2}\bar{\rho}^{2}\|\nabla u\|^{6}_{L^{2}}+M_{0}^{2} M^{3}\bar{\rho}^{3}\|\nabla u\|^{7}_{L^{2}}+M_{0}^{2} M^{4}\bar{\rho}^{4}\|\nabla u\|^{8}_{L^{2}}
	\\
	&\quad+M_{0}^{2}M^{2}\bar{\rho}^{2}\|\nabla u\|^{8}_{L^{2}}+M_{0}^{2} M^{3}\bar{\rho}^{3}\|\nabla u\|^{9}_{L^{2}}+M_{0}^{2} M^{4}\bar{\rho}^{4}\|\nabla u\|^{10}_{L^{2}})ds+\int_{0}^{t}\|\sqrt{\rho}u_{t}\|_{L^{2}}^{2}ds\\
&	\leq\int_{0}^{t}C\mathcal{Q}_{2}(s\|\sqrt{\rho}u_{t}\|_{L^{2}}^{2}+2\underline{\mu}\int_{0}^{s}\tau\int_{\Omega}|\nabla u_{t}|^{2}d\tau)ds+t\|\nabla u\|_{L^{2}}^{2}\int_{0}^{t}C\mathcal{Q}_{1}ds+\int_{0}^{t}\|\sqrt{\rho}u_{t}\|_{L^{2}}^{2}ds,
	\end{aligned}
\end{equation*}
where
\begin{align*}
	\mathcal{Q}_{1}&=M_{0}^{2}M^{2}\bar{\rho}^{2}\|\nabla u\|^{4}_{L^{2}}+M_{0}^{2} M^{3}\bar{\rho}^{3}\|\nabla u\|^{5}_{L^{2}}+M_{0}^{2} M^{4}\bar{\rho}^{4}\|\nabla u\|^{6}_{L^{2}}\\
	&\quad+M_{0}^{2}M^{2}\bar{\rho}^{2}\|\nabla u\|^{6}_{L^{2}}+M_{0}^{2} M^{3}\bar{\rho}^{3}\|\nabla u\|^{7}_{L^{2}}+M_{0}^{2} M^{4}\bar{\rho}^{4}\|\nabla u\|^{8}_{L^{2}},
\end{align*}
and
$$
\mathcal{Q}_{2}= \bar{\rho}M_{0}^{4}\|\nabla u\|^{4}+\bar{\rho}^{3}\|\nabla u\|^{4}+\bar{\rho}\|\nabla u\|^{2}+M_{0}^{2}\bar{\rho}  M^{2}\|\nabla u\|^{2}+
M_{0}^{2}\bar{\rho}  M^{2}\|\nabla u\|^{4}.
$$
Gronwall's inequality yields
\begin{equation}
\begin{aligned}
&\sup_{t\in[0,T]}t\int_{\Omega} \rho u_{t}^{2}+2\underline{\mu}\int_{0}^{T}(t\int_{\Omega}|\nabla u_{t}|^{2})dt
\\&\leq
\int_{0}^{T}\|\sqrt\rho u_{t}\|_{L^{2}}^{2}dt+\sup_{t\in[0,T]}t\|\nabla u\|_{L^{2}}^{2}\int_{0}^{T}C \mathcal{Q}_{1} dt\\&+\int_{0}^{T}(t\|\nabla u \|_{L^{2}}^{2}\int_{0}^{t}C\mathcal{Q}_{1}ds+\int_{0}^{t}\|\sqrt{\rho}u_{t}\|_{L^{2}}^{2}ds)\cdot C\mathcal{Q}_{2}\cdot exp(\int_{0}^{t}C\mathcal{Q}_{2}ds)dt
\\&\leq(\sup_{t\in[0,T]}t\|\nabla u\|_{L^{2}}^{2}\int_{0}^{T}C \mathcal{Q}_{1} dt+\int_{0}^{T}\|\sqrt\rho u_{t}\|_{L^{2}}^{2}dt )\cdot \int_{0}^{T}C\mathcal{Q}_{2}dt\cdot exp(\int_{0}^{T}C\mathcal{Q}_{2}dt),
\\&\leq(\sup_{t\in[0,T]}t\|\nabla u\|_{L^{2}}^{2}\int_{0}^{T}C \mathcal{Q}_{1} dt+\int_{0}^{T}\|\sqrt\rho u_{t}\|_{L^{2}}^{2}dt )\cdot exp(\int_{0}^{T}C\mathcal{Q}_{2}dt)
\end{aligned}
\end{equation}

By using Lemma $\ref{lemma-1}$ and $(\ref{312})$,
one obtains
\begin{equation}
\begin{aligned}
\int_{0}^{T} C\mathcal{Q}_{1}dt&\leq \int_{0}^{T}C(M_{0}^{2}M^{2}\bar{\rho}^{2}\|\nabla u\|^{4}_{L^{2}}+M_{0}^{2} M^{3}\bar{\rho}^{3}\|\nabla u\|^{5}_{L^{2}}+M_{0}^{2} M^{4}\bar{\rho}^{4}\|\nabla u\|^{6}_{L^{2}}\\
&\quad+M_{0}^{2}M^{2}\bar{\rho}^{2}\|\nabla u\|^{6}_{L^{2}}+M_{0}^{2} M^{3}\bar{\rho}^{3}\|\nabla u\|^{7}_{L^{2}}+M_{0}^{2} M^{4}\bar{\rho}^{4}\|\nabla u\|^{8}_{L^{2}})dt\\
&\leq C( {M}^{4}\bar{\rho}^{3}+M^{5}\bar{\rho}^{4}  + {M}^{6}\bar{\rho}^{5})\|u_{0}\|_{L^{2}}^{2},\\	
\end{aligned}
\end{equation}
and
\begin{equation}
	\begin{aligned}
\int_{0}^{T}C \mathcal{Q}_{2}dt&\leq \int_{0}^{T}C( \bar{\rho}M_{0}^{4}\|\nabla u\|^{4}+\bar{\rho}^{3}\|\nabla u\|^{4}+\bar{\rho}\|\nabla u\|^{2}+M_{0}^{2}\bar{\rho}  M^{2}\|\nabla u\|^{2}+
M_{0}^{2}\bar{\rho}  M^{2}\|\nabla u\|^{4})dt\\
&\leq C( {M}^{4}\bar{\rho}^{2}+
\bar{\rho}^{2} + \bar{\rho}^{4})\|u_{0}\|_{L^{2}}^{2}.
\end{aligned}
\end{equation}
Taking ($\ref{331}$) and ($\ref{314}$) into account and using Poincar\'e inequality, we obtain
\begin{equation}
	\begin{aligned}
		&\sup_{t\in[0,T]}t\int_{\Omega} \rho u_{t}^{2}+2\underline{\mu}\int_{0}^{T}(t\int_{\Omega}|\nabla u_{t}|^{2})dt\\
	&\leq \big[C\bar{\rho}( {M}^{4}\bar{\rho}^{3}+M^{5}\bar{\rho}^{4}  + {M}^{6}\bar{\rho}^{5})\|\nabla u_{0}\|_{L^{2}}^{2}\\
&\quad+(2\ln2+1)\frac{\bar{\mu}}{2}\|\nabla u_{0}\|_{L^{2}}\big]\cdot\exp\big[ C( {M}^{4}\bar{\rho}^{2}+
	\bar{\rho}^{2} + \bar{\rho}^{4})\|u_{0}\|_{L^{2}}^{2}\big]\\&
\leq C ({M}^{4}\bar{\rho}^{4}+M^{5}\bar{\rho}^{5}  + {M}^{6}\bar{\rho}^{6}+1)\|\nabla u_{0}\|_{L^{2}}^{2}\cdot\exp\big[ C( {M}^{4}\bar{\rho}^{2}+
\bar{\rho}^{2} + \bar{\rho}^{4})\big].
	\end{aligned}
\end{equation}
 Next we turn to the $t^{2}$-weighted estimates. To do this, multiplying ($\ref{346}$) by t, one has
 \begin{equation}\label{352}
 \begin{aligned}
&\frac{t}{2\underline{\mu}}\frac{d}{dt}(t\int_{\Omega} \rho u_{t}^{2})+t^{2}\int_{\Omega}|\nabla u_{t}|^{2}\\
 &\leq
 Ct^{2}\| \sqrt\rho u_{t}\|_{L^{2}}^{2}( \bar{\rho}M_{0}^{4}\|\nabla u\|_{L^{2}}^{4}+\bar{\rho}^{3}\|\nabla u\|_{L^{2}}^{4}+\bar{\rho}\|\nabla u\|_{L^{2}}^{2}+M_{0}^{2}\bar{\rho}  M^{2}\|\nabla u\|_{L^{2}}^{2}+
M_{0}^{2}\bar{\rho}  M^{2}\|\nabla u\|_{L^{2}}^{4})\\
&\quad+Ct^{2}(M_{0}^{2}M^{2}\bar{\rho}^{2}\|\nabla u\|^{6}_{L^{2}}+M_{0}^{2} M^{3}\bar{\rho}^{3}\|\nabla u\|^{7}_{L^{2}}+M_{0}^{2} M^{4}\bar{\rho}^{4}\|\nabla u\|^{8}_{L^{2}}
+M_{0}^{2}M^{2}\bar{\rho}^{2}\|\nabla u\|^{8}_{L^{2}}\\
&\quad+M_{0}^{2} M^{3}\bar{\rho}^{3}\|\nabla u\|^{9}_{L^{2}}+M_{0}^{2} M^{4}\bar{\rho}^{4}\|\nabla u\|^{10}_{L^{2}})
+\frac{t}{2\underline{\mu}}\|\sqrt\rho u_{t}\|^{2}.
\end{aligned}
\end{equation}
Similar to the method of $\eqref{346}$, adding the term $\frac{t}{2\underline{\mu}}\int_{\Omega}\rho u_{t}^{2}$ into $(\ref{352})$, we deduce that 
\begin{equation}
\begin{aligned}
&\frac{1}{2\underline{\mu}}\frac{d}{dt}[t^{2}\int_{\Omega} \rho u_{t}^{2}]+t^{2}\int_{\Omega}|\nabla u_{t}|^{2}\\
&\leq
Ct^{2}\| \sqrt\rho u_{t}\|_{L^{2}}^{2}( \bar{\rho}M_{0}^{4}\|\nabla u\|_{L^{2}}^{4}+\bar{\rho}^{3}\|\nabla u\|_{L^{2}}^{4}+\bar{\rho}\|\nabla u\|_{L^{2}}^{2}+M_{0}^{2}\bar{\rho}  M^{2}\|\nabla u\|_{L^{2}}^{2}+
M_{0}^{2}\bar{\rho}  M^{2}\|\nabla u\|_{L^{2}}^{4})\\
&\quad+Ct^{2}\|\nabla u\|^{4}_{L^{2}}(M_{0}^{2}M^{2}\bar{\rho}^{2}\|\nabla u\|^{2}_{L^{2}}+M_{0}^{2} M^{3}\bar{\rho}^{3}\|\nabla u\|^{3}_{L^{2}}+M_{0}^{2} M^{4}\bar{\rho}^{4}\|\nabla u\|^{4}_{L^{2}}\\
&\quad
+M_{0}^{2}M^{2}\bar{\rho}^{2}\|\nabla u\|^{4}_{L^{2}}+M_{0}^{2} M^{3}\bar{\rho}^{3}\|\nabla u\|^{5}_{L^{2}}+M_{0}^{2} M^{4}\bar{\rho}^{4}\|\nabla u\|^{6}_{L^{2}})
+\frac{t}{\underline{\mu}}\|\sqrt\rho u_{t}\|_{L^{2}}^{2}.
\end{aligned}
\end{equation}
Integrating with respect to time on $[0,t]$ and using the fact that  $\mu\geq\underline{\mu}>0$, we get
\begin{equation*}
	\begin{aligned}
		&t^{2}\int_{\Omega}\rho u_{t}^{2}+2\underline{\mu}\int_{0}^{t}s^{2}
		\int_{\Omega}|\nabla u_{t}|^{2}ds\\&\leq\int_{0}^{t}Cs^{2}\|\sqrt{\rho}u_{t}\|_{L^{2}}^{2}\mathcal{Q}_{2}ds
		+ \int_{
		0}^{t}Cs^{2}\|\nabla u\|_{L^{2}}^{4}\mathcal{Q}_{3}ds+\int_{0}^{t}2s\|\sqrt{\rho}u_{t}\|_{L^{2}}^{2}ds\\
		&\leq \int_{0}^{t}C\mathcal{Q}_{2}(s^{2}\|\sqrt{\rho}u_{t}\|^{2}_{L^{2}}+2\underline{\mu}\int_{0}^{s}\tau^{2}\int_{\Omega}|\nabla u_{t}|^{2}d\tau )ds+\int_{0}^{t}Cs^{2}\|\nabla u\|_{L^{2}}^{4}\mathcal{Q}_{3}ds+\int_{0}^{t}2s\|\sqrt{\rho}u_{t}\|_{L^{2}}^{2}ds,
	\end{aligned}
\end{equation*}
where \begin{equation*}
	\begin{aligned}
		\mathcal{Q}_{3}&=M_{0}^{2}M^{2}\bar{\rho}^{2}\|\nabla u\|^{2}_{L^{2}}+M_{0}^{2} M^{3}\bar{\rho}^{3}\|\nabla u\|^{3}_{L^{2}}+M_{0}^{2} M^{4}\bar{\rho}^{4}\|\nabla u\|^{4}_{L^{2}}\\
		&\quad
		+M_{0}^{2}M^{2}\bar{\rho}^{2}\|\nabla u\|^{4}_{L^{2}}+M_{0}^{2} M^{3}\bar{\rho}^{3}\|\nabla u\|^{5}_{L^{2}}+M_{0}^{2} M^{4}\bar{\rho}^{4}\|\nabla u\|^{6}_{L^{2}}.
	\end{aligned}
\end{equation*}
Applying Gronwall's inequality, Lemma $\ref{lemma-1}$, ($\ref{312}$) and Lemma $\ref{the5}$, we obtain that
\begin{equation}
\begin{aligned}
&\sup_{t\in[0,T]}t^{2}\int_{\Omega} \rho u_{t}^{2}+2\underline{\mu}\int_{0}^{T}(t^{2}\int_{\Omega}|\nabla u_{t}|^{2})dt\\
&\leq \sup_{t\in[0,T]}t^{2}\|\nabla u\|_{L^{2}}^{4}\int_{0}^{T}C\mathcal{Q}_{3}dt+2\int_{0}^{T}t\|\sqrt{\rho}u_{t}\|_{L^{2}}^{2}dt\\&\quad+(\sup_{t\in[0,T]}t^{2}\|\nabla u\|_{L^{2}}^{4}\int_{0}^{T}C\mathcal{Q}_{3}dt+\int_{0}^{T}t\|\sqrt{\rho}u_{t}\|_{L^{2}}^{2}dt)\cdot\int_{0}^{T}C\mathcal{Q}_{2}dt\cdot \exp(\int_{0}^{T}C\mathcal{Q}_{2}dt)\\&\leq(\sup_{t\in[0,T]}t^{2}\|\nabla u\|_{L^{2}}^{4}\int_{0}^{T}C\mathcal{Q}_{3}dt+\int_{0}^{T}t\|\sqrt{\rho}u_{t}\|_{L^{2}}^{2}dt)\cdot \exp(\int_{0}^{T}C\mathcal{Q}_{2}dt)\\&\leq
C (M^{4}\bar{\rho}^{5}+ M^{5}\bar{\rho}^{6}+ M^{6}\bar{\rho}^{7}+
\bar{\rho})\|\nabla u_{0}\|_{L^{2}}^{2}\cdot\exp\big[ C( {M}^{4}\bar{\rho}^{2}+
\bar{\rho}^{2} + \bar{\rho}^{4})\big],
\end{aligned}
\end{equation}
 where we use \begin{equation*}
 	\int_{0}^{T}C\mathcal{Q}_{3}dt\leq C(M^{4}\bar{\rho}^{3}+ M^{5}\bar{\rho}^{4}+ M^{6}\bar{\rho}^{5}
) \|\nabla u_{0}\|_{L^{2}}^{2}.
 \end{equation*}
Thus, we complete the proof of this lemma.
$\hfill\Box$

We also have the following decay estimates.
\begin{Lemma}\label{lem36}
Assume that $(\rho, u, w, P)$ is the unique strong solution to the system $(\ref{equ11})$-$(\ref{144})$ on $\Omega\times(0,T)$ with the initial data satisfying the assumptions in Theorem \ref{the-11} and (\ref{311})-(\ref{313}).  Then
it holds that
\begin{equation}\label{c}
\sup_{t\in[0,T]}e^{\sigma t}\int_{\Omega}\rho u^{2}(t)+\underline{\mu}\int_{0}^{T}(e^{\sigma t}\int_{\Omega}|\nabla u|^{2})dt\leq  \bar{\rho}\|u_{0}\|_{L^{2}}^{2},
\end{equation}
and
\begin{equation}\label{d}
\int_{0}^{T}(e^{\sigma t}\frac{2}{\underline{\mu}}\int_{\Omega}\rho u_{t}^{2})dt+\sup_{t\in[0,T]}e^{\sigma t}\int_{\Omega}|\nabla u|^{2}\leq C(1+\frac{\bar{\rho}}{\|\rho_{0}\|_{L^{p}}})\|u_{0}\|_{L^{2}}^{2},
\end{equation}
where $\sigma=\frac{\underline{\mu}}{C\|\rho_{0}\|_{L^{p}}}$ for any $p>1$.
\end{Lemma}
{\bf Proof.}
Multiplying the horizontal momentum equations $\eqref{equ11}_{2}$ by $u$ and integrating by parts over $\Omega$, we obtain
\begin{equation}\label{3.511}
\frac{1}{2}\frac{d}{dt}\|\sqrt{\rho}u\|_{L^{2}}^{2}+\int_{\Omega}\mu(\rho)|\nabla u|^{2}=0.
\end{equation}
Moreover, noticing that
\begin{equation}\label{3.51}
\begin{aligned}
\|\sqrt{\rho}u\|_{L^{2}}^{2}&=\int_{\Omega} \rho u^{2}	\leq(\int_{\Omega}\rho^{p})^{\frac{1}{p}}(\int_{\Omega} u^{2q})^{\frac{1}{q}}
\leq\|\rho\|_{L^{p}}\|u\|_{L^{2q}}^{2}
\leq C\|\rho\|_{L^{p}}\|\nabla u\|_{L^{2}}^{2}\\
&\leq C\|\rho_{0}\|_{L^{p}}\|\nabla u\|_{L^{2}}^{2}
\leq\sigma^{-1}\int_{\Omega}\mu(\rho)|\nabla u|^{2},\\
\end{aligned}
\end{equation}
where $\sigma=\frac{\underline{\mu}}{C\|\rho_{0}\|_{L^{p}}}$ for $p>1$.\\
Multiplying  both sides of (\ref{3.511})  and (\ref{3.51}) by the term  $e^{\sigma t}$, respectively, and collecting them, we have
$$\frac{d}{dt}[e^{\sigma t}\|\sqrt{\rho}u\|_{L^{2}}^{2}]+e^{\sigma t}\int_{\Omega} \mu(\rho)|\nabla u|^{2}\leq0.$$
It follows that
\begin{equation}
	\begin{aligned}
		&\sup_{t\in[0,T]}e^{\sigma t}\|\sqrt{\rho}u\|_{L^{2}}^{2}+\underline{\mu}\int_{0}^{T}(e^{\sigma t}\int_{\Omega}|\nabla u|^{2})dt\leq \|\sqrt{\rho_{0}}u_{0}\|_{L^{2}}^{2}\leq  \overline{\rho}\|u_{0}\|_{L^{2}}^{2}.
	\end{aligned}
\end{equation}
Recall that
\begin{equation}\label{360}
2\int_{\Omega}\rho u^{2}_{t}+\frac{d}{dt}\int_{\Omega}\mu(\rho)|\nabla u|^{2}\leq C \mathcal{Q}\|\nabla u\|_{L^{2}}^{2}.
\end{equation}

Multiplying  both sides of (\ref{360}) by the term  $e^{\sigma t}$ and  adding the term $e^{\sigma t} \sigma \int_{\Omega}\mu(\rho)|\nabla u|_{L^{2}}^{2}$ to both sides, we get
\begin{equation}
	\frac{d}{dt}(e^{\sigma t}\int_{\Omega}\mu(\rho)|\nabla u|^{2})+e^{\sigma t}2\int_{\Omega}\rho u_{t}^{2}\leq e^{\sigma t}\|\nabla u\|_{L^{2}}^{2}
C\mathcal{Q}+e^{\sigma t} \sigma \int_{\Omega}\mu(\rho)|\nabla u|^{2}.
\end{equation}
Integrating with respect to time on $[0,t]$ and using the fact that  $\mu\geq\underline{\mu}>0$, we get
\begin{equation*}
	\begin{aligned}
	&	e^{\sigma t}\int_{\Omega}|\nabla u|^{2}+\int_{0}^{t}(\frac{2}{\underline{\mu}}e^{\sigma s}\int_{\Omega}\rho u_{t}^{2})ds\\&\leq \int_{0}^{t}\frac{C\mathcal{Q}}{\underline{\mu}}e^{\sigma s}\|\nabla u\|_{L^{2}}^{2}ds+\int_{0}^{t}e^{\sigma s}\frac{\sigma}{\underline{\mu}}\int_{\Omega}\mu(\rho)|\nabla u|^{2}ds+\int_{\Omega}\frac{\mu(\rho_{0})}{\underline
		{\mu}}|\nabla u_{0}|^{2}\\&\leq
		\int_{0}^{t}\frac{C\mathcal{Q}}{\underline{\mu}}(e^{\sigma s}\|\nabla u\|_{L^{2}}^{2}+\int_{0}^{s}(\frac{2}{\underline{\mu}}e^{\sigma \tau}\int_{\Omega}\rho u_{t}^{2}d\tau))ds\\&\quad+\int_{0}^{t}e^{\sigma s} \frac{\sigma}{\underline{\mu}}\int_{\Omega}\mu(\rho)|\nabla u|^{2}ds+\int_{\Omega}\frac{\mu(\rho_{0})}{\underline{\mu}}|\nabla u_{0}|^{2}
	\end{aligned}
\end{equation*}
By Gronwall inequality in integral form
\begin{equation}
	\begin{aligned}
&\sup_{t\in[0,T]}	e^{\sigma t}\int_{\Omega}|\nabla u|^{2} +	\int_{0}^{T}(e^{\sigma t}\frac{2}{\underline{\mu}}\int_{\Omega}\rho u_{t}^{2})dt\\&\leq \int_{0}^{T}( e^{\sigma
			t}\frac{\sigma}{\underline{\mu}}\int_{\Omega}\mu(\rho)|\nabla u|^{2})dt+\int_{\Omega}\mu(\rho_{0})\frac{|\nabla u_{0}|^{2}}{\underline{\mu}}\\&\quad+[\int_{0}^{T} (e^{\sigma
			t}\frac{\sigma}{\underline{\mu}}\int_{\Omega}\mu(\rho)|\nabla u|^{2})+\int_{\Omega}\mu(\rho_{0})\frac{|\nabla u_{0}|^{2}}{\underline{\mu}}dt]\cdot\int_{0}^{T}\frac{C\mathcal{Q}}{\underline{\mu}}dt\cdot \exp(\int_{0}^{T}\frac{C\mathcal{Q}}{\underline{\mu}}dt)
		\\&\leq(1+\frac{\bar{\rho}}{\|\rho_{0}\|_{L^{p}}})(1+2ln2)\frac{\bar{\mu}}{\underline{\mu}}\|\nabla u_{0}\|_{L^{2}}^{2} .
	\end{aligned}
\end{equation}
Hence, it follows that
\begin{equation}
	\int_{0}^{T}(e^{\sigma t}\frac{2}{\underline{\mu}}\int_{\Omega}\rho u_{t}^{2})dt+\sup_{t\in[0,T]}e^{\sigma t}\int_{\Omega}|\nabla u|^{2}\leq  C(1+\frac{\bar{\rho}}{\|\rho_{0}\|_{L^{p}}})\|\nabla u_{0}\|_{L^{2}}^{2}.
\end{equation}
This completes the proof of Lemma $\ref{lem36}$.$\hfill\Box$

Before stating the $H^{3}$-regularity of the velocity $u$, we  present an estimate for $\int_{0}^{T}\|\nabla u_{t}\|_{L^{2}}^{2}$ as follows.
 \begin{Lemma}\label{LC7}	
 	Assume that $(\rho, u, w, P)$ is the unique strong solution to the system $(\ref{equ11})$-$(\ref{144})$ on $\Omega\times[0,T]$, and the initial data satisfies the assumptions in Theorem \ref{the-11} and (\ref{311})-(\ref{313}). Then, for every $t\in[0,T]$, it holds that
 	\begin{equation}\label{3.64}
 		\|\sqrt{\rho}u_{t}\|_{L^{2}}^{2}+	\underline{\mu}\int_{0}^{t}\|\nabla u_{t}(s)\|_{L^{2}}^{2}ds\leq (\bar{\rho}\|V_{1}\|^{2}_{L^{2}}+CM_{0}^{2}\Theta_{3}\|\nabla u_{0}\|_{L^{2}}^{2})\exp(C\bar{\rho}{\Theta_{3}\|\nabla u_{0}\|^{2}_{L^{2}}}),
 	\end{equation}where $\Theta_{3}=\bar{\rho}M^{2}+\bar{\rho}^{3}M^{2}+\bar{\rho}^{4}M^{3}+\bar{\rho}^{5}M^{4}.$ 
\end{Lemma}

{\bf Proof.}
Multiplying \eqref{335} by $u_{t}$ and integrating by parts over $\Omega$, we have
\begin{equation}\label{3.57}
\begin{aligned}
&\frac{1}{2}\frac{d}{dt}\int_{\Omega}\rho u_{t}^{2}+\int_{\Omega}\mu(\rho)|\nabla u_{t}|^{2}\\&=-\int_{\Omega}\frac{\rho_{t}}{2}u_{t}^{2}-\int_{\Omega}(\rho u)_{t}\partial_{x}u u_{t}+\int_{\Omega}(\rho\int_{0}^{y}\partial_{x}u(s)ds)_{t}\partial_{y}u u_{t}-\frac{1}{2}\int_{\Omega}\mu_{t}(\rho)\frac{d|\nabla u|^{2}}{dt}\\&=
-\int_{\Omega}\frac{\rho_{t}}{2}u_{t}^{2}-\int_{\Omega}\rho_{t}u\partial_{x}u u_{t}-\int_{\Omega}\rho u_{t}\partial_{x}uu_{t}+\int_{\Omega}[\rho_{t}\int_{0}^{y}\partial_{x}u(s)ds\partial_{y}u u_{t}]\\
&\quad+\int_{\Omega}[\rho\int_{0}^{y}\partial_{x}u_{t}(s)ds\partial_{y}u u_{t}]-\int_{\Omega}\mu_{t}(\rho)\nabla u \nabla u_{t}=\Sum_{i=1}^{6}J_{i}.
\end{aligned}
\end{equation}
By direct estimates and Young inequality, we have
\begin{equation}
\begin{aligned}		J_{1}&=\frac{1}{2}\int_{\Omega}(u\partial_{x}\rho-\int_{0}^{y}\partial_{x}u(s)ds\partial_{y}\rho)u_{t}u_{t}
\\&\leq \bar{\rho}^{\frac{1}{2}}\|u\|_{H^{2}}\|\sqrt{\rho}u_{t}\|_{L^{2}}\|\nabla u_{t}\|_{L^{2}}\\&\leq C\bar{\rho} \|u\|_{H^{2}}^{2}\|\sqrt{\rho}u_{t}\|_{L^{2}}^{2}+\frac{~\underline{\mu}~}{2}\|\nabla u_{t}\|_{L^{2}}^{2}.
\end{aligned}
\end{equation}	
It follows from H\"{o}lder's inequality and  Gagliardo-Nirenberg inequality that
\begin{equation}
	\begin{aligned}
		J_{2}&=-\int_{\Omega}\rho_{t}u\partial_{x}u u_{t}=	\|u\|_{H^{2}}\|\partial_{x}\rho\|_{L^{\infty}}\|\partial_{x}u\|_{L^{2}}\|u_{t}\|_{L_{x}^{\infty}L_{y}^{2}}\|u\|_{L_{x}^{2}L_{y}^{\infty}}\\
		&\quad+\|u\|_{H^{2}}\|\partial_{x}\rho\|_{L^{\infty}}\|\partial_{x}u\|_{L^{2}}\|u_{t}\|_{L_{x}^{\infty}L_{y}^{2}}\|\int_{0}^{y}\partial_{x}u(s)ds\|_{L_{x}^{2}L_{y}^{\infty}}\\&\leq C\|\nabla \rho\|_{L^{\infty}}^{2} \|u\|_{H^{2}}^{2}\|\nabla u\|_{L^{2}}^{4}+\frac{~\underline{\mu}~}{2}\|\nabla u_{t}\|_{L^{2}}^{2}.
	\end{aligned}
\end{equation}
Similarly,
\begin{equation}
	\begin{aligned}
		J_{3}&=-\int_{\Omega}	\rho u_{t}\partial_{x}u u_{t}\\
		&\leq\sqrt{\bar{\rho}}\|\sqrt{\rho}u_{t}\|_{L^{2}}\|u_{t}\|_{L_{x}^{\infty}L_{y}^{2}}\|\partial_{x}u\|_{L_{x}^{2}L_{y}^{\infty}}\\
		&\leq C \bar{\rho}\|\sqrt{\rho}u_{t}\|_{L^{2}}^{2}\|u\|_{H^{2}}^{2}+\frac{~\underline{\mu}~}{2}\|\nabla u_{t}\|_{L^{2}}^{2}.
	\end{aligned}
\end{equation}
Applying anisotropic Sobolev embedding theorem  and Young inequality, one has
\begin{equation}
\begin{aligned}	J_{4}&=\int_{\Omega}[(u\partial_{x}\rho-\int_{0}^{y}\partial_{x}u(s)ds\partial_{y}\rho)\int_{0}^{y}\partial_{x}u(s)ds \partial_{y}u u_{t}]\\&\leq\|u\|_{L^{\infty}}\|\partial_{x}\rho\|_{L^{\infty}}\|\int_{0}^{y}\partial_{x}u(s)ds\|_{L^{2}_{x}L^{\infty}_{y}}\|\partial_{y}u\|_{L^{2}}\|u_{t}\|_{L^{\infty}_{x}L^{2}_{y}}\\		&\quad+\|\int_{0}^{y}\partial_{x}u(s)ds\|_{L^{\infty}}\|\partial_{y}\rho\|_{L^{\infty}}\|\int_{0}^{y}\partial_{x}u(s)ds\|_{L^{2}_{x}L^{\infty}_{y}}\|\partial_{y}u\|_{L^{2}}\|u_{t}\|_{L^{\infty}_{x}L^{2}_{y}}\\
&\leq C\|\nabla \rho\|_{L^{\infty}}^{2}\|u\|_{H^{2}}^{2}\|\nabla u\|_{L^{2}}^{4}+\frac{~\underline{\mu}~}{2}\|\nabla u_{t}\|_{L^{2}}^{2},
\end{aligned}
\end{equation}

\begin{equation}
\begin{aligned}	J_{5}
&=\int_{\Omega}[\rho\int_{0}^{y}\partial_{x}u_{t}(s)ds\partial_{y}uu_{t}]\\&\leq\|\int_{0}^{y}\partial_{x}u_{t}(s)ds\|_{L_{x}^{2}L_{y}^{\infty}}\|\partial_{y}u
\|_{L_{y}^{2}L_{x}^{\infty}}\|\rho u_{t}\|_{L^{2}}\\
&\leq C\bar{\rho}\|\sqrt{\rho}u_{t}\|_{L^{2}}^{2}\|u\|_{H^{2}}^{2}+\frac{~\underline{\mu}~}{2}\|\nabla u_{t}\|_{L^{2}}^{2},\\
\end{aligned}
\end{equation}
and
\begin{equation}
	\begin{aligned}
		J_{6}=\int_{\Omega}\mu_{t}(\rho)\nabla u \nabla u_{t}&\leq C\|\nabla \rho\|_{L^{\infty}}^{2}\|u\|_{H^{2}}^{2}\|\nabla u\|_{L^{2}}^{4}+\frac{~\underline{\mu}~}{2}\|\nabla u_{t}\|_{L^{2}}^{2}.
	\end{aligned}	
\end{equation}
Substituting the above estimates $J_{1}$-$J_{6}$ into $(3.57)$  we obtain
\begin{equation}
\frac{d}{dt}\int_{\Omega}\rho u_{t}^{2}+\underline{\mu}\int_{\Omega}|\nabla u_{t}|^{2}\leq C(\bar{\rho}\|u\|_{H^{2}}^{2}\|\sqrt{\rho}u_{t}\|_{L^{2}}^{2}+M_{0}^{2}\|\nabla u\|_{L^{2}}^{4}\|u\|_{H^{2}}^{2}).
\end{equation}
Then, integrating it over $(0,t)\subset(0,T)$ and using Lemma \ref{lemma-1}-\ref{lem-33}, one has
\begin{equation}\label{399}
\begin{aligned}
&	\|\sqrt{\rho}u_{t}\|_{L^{2}}^{2}+\underline{\mu}\int_{0}^{t}\int_{\Omega}|\nabla u_{t}|^{2}ds\\
&\leq [	\|\sqrt{\rho}u_{t}(0)\|_{L^{2}}^{2}+C\int_{0}^{T}(M_{0}^{2}\|\nabla u\|_{L^{2}}^{4}\|u\|_{H^{2}}^{2})ds]\exp(C\bar{\rho}\int_{0}^{T}\|u\|^{2}_{H^{2}}ds)\\
&\leq [\bar{\rho}\|V_{1}\|_{L^{2}}^{2}+CM_{0}^{2}\Theta_{3}\|\nabla u_{0}\|_{L^{2}}^{2}]\exp(C\bar{\rho}{\Theta_{3}\|\nabla u_{0}\|^{2}_{L^{2}}}),
\end{aligned}
\end{equation}
where $\Theta_{3}=\bar{\rho}M^{2}+\bar{\rho}^{3}M^{2}+\bar{\rho}^{4}M^{3}+\bar{\rho}^{5}M^{4}$.
Then, we complete the proof of this lemma
.$\hfill\Box$\\
\begin{Remark}
	When we deal with the term of 
~$\|\sqrt{\rho}u_{t}(0)\|_{L^{2}}$, we first need to regularize the initial data ($\rho_{0}, u_{0}, P_{0}$). That is, for each $\delta\in(0,1)$, that we choose $\rho^{\delta}_{0} \in C^{2+\alpha}$ such that $0<\delta\leq\rho^{\delta}_{0}\leq \rho_{0}+1$, $\rho_{0}^{\delta}\rightarrow\rho_{0}$ in $W^{2,2}$ as $\delta\rightarrow0$, and denote by $( u_{0}^{\delta},P_{0}^{\delta})\in H^{1}(\Omega)\times L^{2}(\Omega)$ a solution to the problem \begin{equation}\label{4.10}
		-\rho^{\delta}_{0}u^{\delta}_{0}\partial_{x}u^{\delta}_{0}+\rho^{\delta}_{0}\int_{0}^{y}\partial_{x}u^{\delta}_{0}(s)ds\partial_{y}u^{\delta}_{0}-\partial_{x}P^{\delta}_{0}+\partial_{x}(\mu(\rho^{\delta}_{0})\partial_{x}u^{\delta}_{0})+\partial_{y}(\mu(\rho^{\delta}_{0})\partial_{y}u^{\delta}_{0}):=\rho^{\delta}_{0}V_{1}.
	\end{equation}In fact, with such regularized initial data and the compatibility condition (\ref{4.10}), we can get corresponding solutions $(\rho^{\delta},u^{\delta}, P^{\delta})$ satisfying the estimates in Lemma $\ref{lemma-1}$-$\ref{pro-41}$ and these estimates are independent of $\delta$. Then we use the standard compactness arguments to get the desired result.
\end{Remark}
With the help of Lemmas \ref{the5}-\ref{LC7}, we can now estimate
$\|u\|_{L^{1}(0,T ;H^{3}(\Omega))}$.
\begin{Lemma}\label{the8}
	Assume that $(\rho, u, w, P)$ is the unique strong solution to the system $(\ref{equ11})$-$(\ref{144})$ on $\Omega\times(0,T)$ with the initial data satisfying the assumptions in Theorem \ref{the-11} and (\ref{311})-(\ref{313}). Then it holds that
\begin{equation}\label{j}
	\begin{aligned}
		\int_{0}^{T}\|u\|_{H^{3}}dt&\leq\tilde{C}\Theta_{5} (\bar{\rho}C+M_{0}^{2}\Theta_{3}\|\nabla u_{0}\|_{L^{2}}^{2})\exp(C\bar{\rho}{\Theta_{3}\|\nabla u_{0}\|^{2}_{L^{2}}})\sqrt{T_{0}}\\
		&\quad+\tilde{C}\Theta_{5}(\bar{\rho}\|\nabla u_{0}\|_{L^{2}}^{2} \cdot\Theta_{1}\cdot\exp(\Theta_{2}))^{\frac{1}{2}}\frac{1}{\sqrt{T_{0}}}+\tilde{C}M^{5}\bar{\rho}^{2}\|\nabla u_{0}\|_{L^{2}}^{2}+\tilde{C}\Theta_{6}\bar{\rho}\|\nabla u_{0}\|_{L^{2}}^{2},
	\end{aligned}
\end{equation}
where $\Theta_{5}=M^{5}\bar{\rho}+M^{4}$, $\Theta_{6}=M^{5}\bar{\rho}+M^{5}\bar{\rho}^{3}+M^{\frac{11}{2}}\bar{\rho}^{\frac{3}{2}}+M^{6}\bar{\rho}^{2}+M^{6}\bar{\rho}^{4}+M^{7}\bar{\rho}^{5}$,
and $T_{0}<T$ which is to be determined later.
\end{Lemma}

{\bf Proof.}
It follows from the regularity theory for the density-dependent hydrostatic Stokes system (see Lemma \ref{B2}) that
\begin{equation}
\begin{aligned}
\|u\|_{H^{3}}&
\leq\tilde{C}M^{3}\|-\rho u_{t}+\rho u\partial_{x}u+\rho\int_{0}^{y}\partial_{x}u(s)ds\partial_{y}u\|_{H^{1}}\triangleq L_{1}+L_{2},
\end{aligned}
\end{equation}
where
$$L_{1}=\tilde{C}M^{3} (\|\sqrt\rho\sqrt\rho u_{t}\|_{L^{2}}+\|\rho u\partial_{x}u\|_{L^{2}}+\|\rho \int_{0}^{y}\partial_{x}u(s)ds\partial_{y}u\|_{L^{2}}),$$
and
\begin{equation}
	\begin{aligned}
	L_{2}=&\tilde{C}M^{3}(\|\nabla \rho u_{t}\|_{L^{2}}+\|\rho \nabla u_{t}\|_{L^{2}}+\|\nabla \rho u\partial_{x}u\|_{L^{2}} +\|\rho \nabla u\partial_{x}u\|_{L^{2}} +\|\rho u\nabla\partial_{x}u\|_{L^{2}}
	\\&+\|\nabla\rho \int_{0}^{y}\partial_{x}u(s)ds\partial_{y}u\|_{L^{2}} +\|\rho \partial_{x}u\partial_{y}u\|_{L^{2}}+\| \rho \int_{0}^{y}\partial_{xx}u(s)ds\partial_{y}u\|_{L^{2}}
\\&	+\|\rho\int_{0}^{y}\partial_{x}u(s)ds\nabla\partial_{y}u\|_{L^{2}}).
	\end{aligned}
\end{equation}
Similar to the proof of $\eqref{CC}$ in Lemma \ref{lem-33}, we obtain
\begin{equation}\label{367}
\begin{aligned}
L_{1}
&\leq\tilde{C}M^{3}(\bar{\rho}^{1/2}\| \sqrt\rho u_{t}\|_{L^{2}}+ \bar{\rho}\|\nabla u\|^{2}_{L^{2}}+  \bar{\rho} \|\nabla u\|_{L^{2}}^{\frac{5}{3}}\|\nabla ^{2}u\|_{L^{2}}^{\frac{1}{3}}+\bar{\rho} \|\nabla u\|_{L^{2}}^{\frac{3}{2}}\|\nabla ^{2}u\|_{L^{2}}^{\frac{1}{2}})\\
&\leq \frac{1}{2}\|\nabla u\|_{H^{2}}+ \tilde{C}M^{3}\bar{\rho}^{1/2}\| \sqrt\rho u_{t}\|_{L^{2}}+ \tilde{C}\bar{\rho}M^{3}\|\nabla u\|^{2}_{L^{2}}+  \tilde{C}\bar{\rho}^{\frac{3}{2}} M^{\frac{9}{2}}\|\nabla u\|_{L^{2}}^{\frac{5}{2}}+\tilde{C}\bar{\rho}^{2} M^{6} \|\nabla u\|_{L^{2}}^{3}.
	\end{aligned}
\end{equation}
\\
We rewrite $L_{2}$ as follow
\begin{equation}\label{C50}
\begin{aligned}
L_{2}&\leq\tilde{C}M^{3}\|\nabla \rho u_{t}\|_{L^{2}} +\tilde{C}M^{3}\|\rho\nabla u_{t}\|_{L^{2}} +\tilde{C}M^{3}\|\nabla \rho u\partial_{x}u\|_{L^{2}} +\tilde{C}M^{3}\| \rho \nabla u\nabla u\|_{L^{2}} \\
&\quad+\tilde{C}M^{3}\|\rho u \nabla^{2}u\|_{L^{2}} +\tilde{C}M^{3}\|\nabla\rho \int_{0}^{y}\partial_{x}u(s)ds \partial_{y}u\|_{L^{2}} +\tilde{C}M^{3}\|\rho\int_{0}^{y}\partial_{xx}u(s)ds\partial_{y}u\|_{L^{2}} \\
&\quad+\tilde{C}M^{3}\|\rho\int_{0}^{y}\partial_{x}u(s)ds\nabla \partial_{y}u \|_{L^{2}})=\sum_{i=1}^{8}J_{2i}.
\end{aligned}
\end{equation}
By direct estimates, one has
\begin{equation}
	L_{21}+L_{22}\leq\tilde{C}M^{4} \|\nabla u_{t}\|_{L^{2}}	+\tilde{C}M^{3}\bar{\rho}\|\nabla u_{t}\|_{L^{2}}.
\end{equation}
Applying H\"{o}lder's inequality, $\eqref{318}$ and Gagliardo-Nirenberg inequality, we deduce that
\begin{equation}
\begin{aligned}
L_{23}
&\leq\tilde{C}M^{3} M_{0}\|u\|_{H^{2}}\|\partial_{x}u\|_{L^{2}}\\
&\leq \tilde{C}M^{3}M_{0}\|\nabla u\|_{L^{2}}( \bar{\rho}^{1/2}M\| \sqrt\rho u_{t}\|_{L^{2}}+ \bar{\rho}M\|\nabla u\|^{2}_{L^{2}}+  \bar{\rho}^{\frac{3}{2}}M^{\frac{3}{2}}\|\nabla u\|_{L^{2}}^{\frac{5}{2}}+\bar{\rho}^{2} M^{2} \|\nabla u\|_{L^{2}}^{3})\\
&\leq\tilde{C}M^{5}\bar{\rho}^{\frac{1}{2}}\|\sqrt{\rho }u_{t}\|_{L^{2}}\|\nabla u\|_{L^{2}}+\tilde{C}M^{5}\bar{\rho}\|\nabla u\|_{L^{2}}^{3}+\tilde{C}M^{\frac{11}{2}}\bar{\rho}^{\frac{3}{2}}\|\nabla u\|_{L^{2}}^{\frac{7}{2}}+\tilde{C}M^{6}\bar{\rho}^{2}\|\nabla u\|_{L^{2}}^{4},
\end{aligned}
\end{equation}
\begin{equation}
\begin{aligned}
L_{24}
&\leq\tilde{C}M^{3}\bar{\rho}\|\nabla u\|_{L_{y}^{2}L_{x}^{\infty}}\|\nabla u\|_{L_{x}^{2}L_{y}^{\infty}}
\leq\tilde{C}M^{3}\bar{\rho}\|u\|_{H^{2}}^{2}\\&\leq \tilde{C}M^{3}\bar{\rho}( \bar{\rho}^{1/2}M\| \sqrt\rho u_{t}\|_{L^{2}}+ \bar{\rho}M\|\nabla u\|^{2}_{L^{2}}+  \bar{\rho}^{\frac{3}{2}}M^{\frac{3}{2}}\|\nabla u\|_{L^{2}}^{\frac{5}{2}}+\bar{\rho}^{2} M^{2} \|\nabla u\|_{L^{2}}^{3})^{2}\\
&\leq \tilde{C}M^{3}\bar{\rho}( \bar{\rho}M^{2}\| \sqrt\rho u_{t}\|^{2}_{L^{2}}+\bar{\rho}^{2}M^{2}\|\nabla u\|^{4}_{L^{2}}+  \bar{\rho}^{3}M^{3}\|\nabla u\|_{L^{2}}^{5}+\bar{\rho}^{4} M^{4} \|\nabla u\|_{L^{2}}^{6}\\
&\leq\tilde{C}M^{5}\bar{\rho}^{2}\|\sqrt{\rho }u_{t}\|_{L^{2}}^{2}+\tilde{C}M^{5}\bar{\rho}^{3}\|\nabla u\|_{L^{2}}^{4}+\tilde{C}M^{6}\bar{\rho}^{4}\|\nabla u\|_{L^{2}}^{5}+\tilde{C}M^{7}\bar{\rho}^{5}\|\nabla u\|_{L^{2}}^{6},
\end{aligned}
 \end{equation}
and
\begin{equation}
	L_{25}\leq\tilde{C}M^{3}\bar{\rho}\|u\|_{H^{2}}^{2}.
\end{equation}
By Sobolev and Poincar\'{e} inequalities, We estimate $I_{26}$ as follows
\begin{equation}
\begin{aligned}
	L_{26}
&\leq\tilde{C}M^{3}\|\nabla \rho\|_{L^{\infty}}\|\int_{0}^{y}\partial_{x}u(s)ds \|_{L_{y}^{\infty}L_{x}^{2}}\|\partial_{y}u \|_{L_{x}^{\infty}L_{y}^{2}}\leq\tilde{C}M^{3} M_{0}\|\nabla u\|_{L^{2}}\|u\|_{H^{2}}.
\end{aligned}
\end{equation}
Similarly, one has
\begin{equation}
	L_{27}+L_{28}\leq\tilde{C}M^{3}\bar{\rho}\|u\|_{H^{2}}^{2}.
\end{equation}
Substituting  the above estimates $L_{21}$-$L_{28}$ into $\eqref{C50}$, we have
\begin{equation}\label{376}
\begin{aligned}
L_{2}\leq&\tilde{C}M^{4} \|\nabla u_{t}\|_{L^{2}}	+\tilde{C}M^{3}\bar{\rho}\|\nabla u_{t}\|_{L^{2}}+\tilde{C}M^{5}\bar{\rho}^{2}\|\sqrt{\rho }u_{t}\|_{L^{2}}^{2}\\&   
+\tilde{C}M^{5}\bar{\rho}^{\frac{1}{2}}\|\sqrt{\rho }u_{t}\|_{L^{2}}\|\nabla u\|_{L^{2}}+\tilde{C}M^{5}\bar{\rho}\|\nabla u\|_{L^{2}}^{3}+\tilde{C}M^{\frac{11}{2}}\bar{\rho}^{\frac{3}{2}}\|\nabla u\|_{L^{2}}^{\frac{7}{2}}\\
&+\tilde{C}M^{6}\bar{\rho}^{2}\|\nabla u\|_{L^{2}}^{4}+\tilde{C}M^{5}\bar{\rho}^{3}\|\nabla u\|_{L^{2}}^{4}+\tilde{C}M^{6}\bar{\rho}^{4}\|\nabla u\|_{L^{2}}^{5}+\tilde{C}M^{7}\bar{\rho}^{5}\|\nabla u\|_{L^{2}}^{6}.
	\end{aligned}
\end{equation}
Combining (\ref{367}) and (\ref{376}), one has
\begin{equation}\label{395}
\begin{aligned}
\int_{0}^{T}\|u\|_{H^{3}} dt\leq& \int_{0}^{T}[ \tilde{C}M^{3}\bar{\rho}^{1/2}\| \sqrt\rho u_{t}\|_{L^{2}}
+\tilde{C}M^{4} \|\nabla u_{t}\|_{L^{2}}	+\tilde{C}M^{3}\bar{\rho}\|\nabla u_{t}\|_{L^{2}}\\&+\tilde{C}M^{5}\bar{\rho}^{\frac{1}{2}}\|\sqrt{\rho }u_{t}\|_{L^{2}}\|\nabla u\|_{L^{2}}+\tilde{C}M^{5}\bar{\rho}^{2}\|\sqrt{\rho }u_{t}\|_{L^{2}}^{2}+ \tilde{C}\bar{\rho}M^{3}\|\nabla u\|^{2}_{L^{2}} \\&	+ \tilde{C}\bar{\rho}^{\frac{3}{2}} M^{\frac{9}{2}}\|\nabla u\|_{L^{2}}^{\frac{5}{2}}+\tilde{C}\bar{\rho}^{2} M^{6} \|\nabla u\|_{L^{2}}^{3}
+\tilde{C}M^{5}\bar{\rho}\|\nabla u\|_{L^{2}}^{3}+\tilde{C}M^{\frac{11}{2}}\bar{\rho}^{\frac{3}{2}}\|\nabla u\|_{L^{2}}^{\frac{7}{2}}\\
&+\tilde{C}M^{6}\bar{\rho}^{2}\|\nabla u\|_{L^{2}}^{4}+\tilde{C}M^{5}\bar{\rho}^{3}\|\nabla u\|_{L^{2}}^{4}+\tilde{C}M^{6}\bar{\rho}^{4}\|\nabla u\|_{L^{2}}^{5}+\tilde{C}M^{7}\bar{\rho}^{5}\|\nabla u\|_{L^{2}}^{6}]dt\\
\leq&
\tilde{C}\Theta_{5}\int_{0}^{T}\|\nabla u_{t}\|_{L^{2}}dt
+\tilde{C}M^{5}(\bar{\rho}^{2}+1)\int_{0}^{T}\|\sqrt{\rho }u_{t}\|_{L^{2}}^{2}dt+\tilde{C}
\Theta_{6}\int_{0}^{T}\|\nabla u\|_{L^{2}}^{2}dt,
\end{aligned}
\end{equation}
where $\Theta_{5}=M^{3}\bar{\rho}+M^{4}+M^{3}\bar{\rho}^{\frac{1}{2}}$ and $\Theta_{6}=M^{5}\bar{\rho}+M^{5}\bar{\rho}^{3}+M^{\frac{11}{2}}\bar{\rho}^{\frac{3}{2}}+M^{6}\bar{\rho}^{2}+M^{6}\bar{\rho}^{4}+M^{7}\bar{\rho}^{5}$.Then, in order to obtain the estimate of $\|u\|_{L^{1}(0,T ;H^{3}(\Omega))}$, we have to deal with the term of
$$\int_{0}^{T} \| \nabla u_{t}(t)\|_{L^{2}}dt.$$\\
Applying H\"{o}lder's inequality, (\ref{3.34}) and (\ref{3.64}), one has
\begin{equation}\label{393}
\begin{aligned}	
\int_{0}^{T}\|\nabla u_{t}\|_{L^{2}}dt
&\leq(\int_{0}^{T_{0}}\|\nabla u_{t}\|_{L^{2}}^{2}dt)^{\frac{1}{2}}\sqrt{T_{0}}+(\int_{T_{0}}^{T}t^{2}\|\nabla u_{t}\|_{L^{2}}^{2}dt)^{\frac{1}{2}} (\int_{T_{0}}^{T}t^{-2}dt)^{\frac{1}{2}}\\
&\leq C[\bar{\rho}+M_{0}^{2}\Theta_{3}\|\nabla u_{0}\|^{2}_{L^{2}}]\exp(C\bar{\rho}\Theta_{3}\|\nabla u_{0}\|^{2}_{L^{2}})\sqrt{T_{0}}\\&\quad+C[\bar{\rho}\|\nabla u_{0}\|_{L^{2}}^{2} \cdot\Theta_{1}\cdot\exp(\Theta_{2})]^{\frac{1}{2}}\frac{1}{\sqrt{T_{0}}},
\end{aligned}
\end{equation}
where $T_{0} < T$ is a small number which is determined later.

Thus, substituting $\eqref{393}$ into \eqref{395} and using Poincar\'{e} inequalities, we obtain that
\begin{equation}\label{391}
\begin{aligned}
\int_{0}^{T}\|u\|_{H^{3}}dt
\leq&
\tilde{C}\Theta_{5}[\bar{\rho}+M_{0}^{2}\Theta_{3}\|\nabla u_{0}\|^{2}_{L^{2}}]\exp(C\bar{\rho}\Theta_{3}\|\nabla u_{0}\|^{2}_{L^{2}})\sqrt{T_{0}}\\&+\tilde{C}\Theta_{5}[\bar{\rho}\|\nabla u_{0}\|_{L^{2}}^{2} \cdot\Theta_{1}\cdot\exp(\Theta_{2})]^{\frac{1}{2}}\frac{1}{\sqrt{T_{0}}}\\&
+\tilde{C}M^{5}(\bar{\rho}^{2}+1)\|\nabla u_{0}\|_{L^{2}}^{2}+\tilde{C}
\Theta_{6}\bar{\rho}\|\nabla u_{0}\|_{L^{2}}^{2}.
\end{aligned}
\end{equation}
Here the proof of this lemma is completed.
$\hfill\Box$

In view of (\ref{391}), we are now in a position to prove that $\|\nabla \rho(t)\|_{L^{\infty
}}+\|\nabla^{2} \rho(t)\|_{L^{2}}\leq 4M_{0}$,  
provided the initial data $\|\nabla u_{0}\|_{L^{2}}$ and the time $T_{0}$ are suitably small in some sense.

 \begin{Lemma}\label{pro-41}
Assume that $(\rho, u, w, P)$ is the unique strong solution to the system $(\ref{equ11})$-$(\ref{144})$ on $\Omega\times[0,T]$, and the initial data satisfies the assumptions in Theorem \ref{the-11} and (\ref{311})-(\ref{313}). Moreover, suppose that there exists a positive number $\epsilon_{0}$, depending on $\Omega$, $\sup_{\Omega}\rho_{0}=\overline{\rho}$,  $\|\mu\|_{C^{2}}$, $\|\nabla\rho_{0}\|_{L^{\infty}}$, $\| \nabla^{2} \rho_{0}\|_{L^{2}}$, such that
$$ \|\nabla u_{0}\|_{L^{2}}\leq \epsilon_{0}.$$
Then it holds that
\begin{equation}\label{100}
	\sup_{t\in[0,T]}\|\nabla \rho(t)\|_{L^{\infty}}+\sup_{t\in[0,T]}\|\nabla^{2} \rho(t)\|_{L^{2}}\leq 4M_{0}.
\end{equation}
\end{Lemma}

{\bf Proof.}  Recall that
\begin{equation}\label{42}
	\|\nabla \rho\|_{L^{\infty}}\leq \|\nabla \rho_{0}\|_{L^{\infty}}\exp(\int_{0}^{t}\|u(s)\|_{H^{3}}ds),
\end{equation}
and
\begin{equation}\label{43}
	\|\nabla ^{2}\rho\|_{L^{2}}\leq (\|\nabla \rho_{0}\|_{L^{\infty}}+\|\nabla ^{2}\rho_{0}\|_{L^{2}})\exp(\int_{0}^{t}\|u(s)\|_{H^{3}}ds).
\end{equation}Now, we can choose some small positive constant $\epsilon_{0}$ and then time $T_{0}<T$, satisfying
\begin{equation}\label{jj}
\begin{aligned}
	&\tilde{C}\Theta_{5} [\bar{\rho}+M_{0}^{2}\Theta_{3}\epsilon_{0}^{2}]\exp(C\bar{\rho}\Theta_{3}\epsilon_{0}^{2})\sqrt{T_{0}}\\
	&+\tilde{C}\Theta_{5}[\bar{\rho}\epsilon_{0}^{2} \cdot\Theta_{1}\cdot\exp(\Theta_{2})]^{\frac{1}{2}}\frac{1}{\sqrt{T_{0}}}+\tilde{C}M^{5}\bar{\rho}^{2}\epsilon_{0}^{2}+\tilde{C}\Theta_{6}\bar{\rho}\epsilon_{0}^{2}\leq\ln2
\end{aligned}
\end{equation}
Thus, if $\|\nabla u_{0}\|_{L^{2}}\leq \epsilon_{0}$, it follows from Lemma \ref{the8} that
$$\|\nabla\rho\|_{L^{\infty}}+ \|\nabla^{2}\rho\|_{L^{2}}\leq 2(2\|\rho_{0}\|_{L^{\infty}}+\|\nabla^{2}\rho_{0}\|_{L^{2}})\leq 4 (\|\nabla \rho_{0}\|_{L^{\infty}}+\|\nabla^{2}\rho_{0}\|_{L^{2}})=4M_{0}.$$
Here we complete the proof of Lemma \ref{pro-41}.
$\hfill\Box$

\subsection{ Proof of Theorem \ref{the-2}}
With the above uniform estimates in hand, we are  ready to prove Theorem \ref{the-2}.

{\bf Proof of Theorem \ref{the-2}.} According to Lemma \ref{coro1}, there exists a $T_{*}>0$ such that the primitive equations have a unique local strong solution $(\rho, u, w, P)$ on $[0,T_{*}]$.  Next, we apply the above a priori estimates in Section 3.1 to extend the local solution to a global one.

Since $\|\nabla \rho_{0}\|_{L^{\infty}}+ \|\nabla^{2} \rho_{0}\|_{L^{2}}=M_{0}<8M_{0}$, and thanks to the continuity of $\|\nabla\rho(t)\|_{L^{\infty}}+\|\nabla^{2} \rho(t)\|_{L^{2}}$  and $\|\nabla u(t)\|_{L^{2}}$, there exists a $T_{1} \in (0,T_{*})$ such that
$$\sup_{0\leq t\leq T_{1}}\|\nabla\rho(t)\|_{L^{\infty}}+\|\nabla^{2} \rho(t)\|_{L^{2}}\leq 8M_{0},$$
and 
$$\sup_{0\leq t\leq T_{1}}\|\nabla u(t)\|_{L^{2}}^{2}\leq 4\frac{\bar{\mu}}{\underline{\mu}}\|\nabla u_{0}\|_{L^{2}}^{2}.$$
Let
$$ T^{*}=sup\Bigl\{T|(\rho, u, w, P) ~is~a~ strong~ solution~ to~(1.1 )~ on~ [0,T]\Bigr\},$$
and
\[T_{1}^{*}=sup\left\{T\Bigg|
\begin{aligned}
&(\rho, u, w,P)\text{~is~ a~ strong~ solution~ to~ (1.1 )~ on~} [0,T],\sup_{0\leq t\leq T}(\|\nabla\rho(t)\|_{L^{\infty}}\\	
&+\|\nabla^{2} \rho(t)\|_{L^{2}})\leq 8M_{0}, \text{~and ~}\sup_{0\leq t\leq T}\|\nabla u(t)\|_{L^{2}}^{2}\leq (2ln2+1)\frac{\bar{\mu}}{\underline{\mu}}\|\nabla u_{0}\|_{L^{2}}^{2}\
\end{aligned}
\right\}.\]
Then $T^{*}_{1}\geq T_{1}>0.$ Recalling Lemma \ref{prop-1} and Lemma \ref{pro-41}, it is easy to verify that $$T^{*}=T^{*}_{1},$$ provided that $\|\nabla u_{0}\|_{L^{2}}\leq \eta_{0}=\min\{\gamma_{0},\epsilon_{0}\}$ as assumed.

We claim that $T^{*}=\infty.$ Otherwise assume that $T^{*}<\infty.$ By virtue of Lemma \ref{prop-1} and Lemma \ref{pro-41}, for every $t\in[0,T^{*})$, it holds that
$$\|\nabla\rho(t)\|_{L^{\infty}}+\|\nabla^{2} \rho\|_{L^{2}}\leq 4M_{0} ~\text{and}~\|\nabla u(t)\|_{L^{2}}^{2}\leq(2ln2+1)\frac{\bar{\mu}}{\underline{\mu}}\|\nabla u_{0}\|_{L^{2}}^{2},$$
which contradicts the blow-up criterion (\ref{16}) in Theorem \ref{the-11}.
Hence we complete the proof of Theorem \ref{the-2}.
\  \  \

 {\bf Acknowledgements.} Q. Jiu is partially
supported by National Natural Sciences Foundation of China (No. 11931010).
F. Wang is partially supported by the National Natural Science Foundation of China (No.12201028), Chinese Postdoctoral Science Foundation (No.2022M720383)

\end{document}